\documentclass[11pt]{article}
\usepackage{anysize}
\marginsize{2.0cm}{2.0cm}{2.0 cm}{2.0 cm}
\usepackage{setspace}

\usepackage{latexsym,amsfonts,amsmath}
\usepackage{graphicx}
\usepackage{color}
\usepackage{epstopdf}
\usepackage[titletoc]{appendix}

\def\limsup{\mathop{\overline{\lim}}}

\newtheorem{condition}{Condition}[section]{\bfseries}{\itshape}

{\bfseries}{\itshape}

\newtheorem{theorem}{Theorem}[section]{\bfseries}{\itshape}

\newtheorem{corollary}{Corollary}[section]{\bfseries}{\itshape}

\newtheorem{proposition}{Proposition}[section]{\bfseries}{\itshape}

\newtheorem{example}{Example}[section]{\bfseries}{\itshape}

\newtheorem{lemma}{Lemma}[section]{\bfseries}{\itshape}

\newtheorem{remark}{Remark}[section]{\bfseries}{\itshape}

\newtheorem{definition}{Definition}[section]{\bfseries}{\itshape}

{\bfseries}{\itshape}

\begin{document}

\title{Constrained total undiscounted continuous-time Markov decision processes}
\author{Xianping Guo\thanks{School of Mathematics and Computational Science, Sun Yat-Sen University,
Guangzhou, P.R. China. E-mail: mcsgxp@mail.sysu.edu.cn.} ~and Yi
Zhang \thanks{Department of Mathematical Sciences, University of
Liverpool, Liverpool, L69 7ZL, U.K.. E-mail: yi.zhang@liv.ac.uk.}}

\maketitle

\par\noindent{\bf Abstract:}
The present paper considers the constrained optimal control
problem with total undiscounted criteria for a continuous-time
Markov decision process (CTMDP) in Borel state and action spaces.
Under the standard compactness and continuity conditions, we show
the existence of an optimal stationary policy out of the class of
general nonstationary ones. In the process, we justify the
reduction of the CTMDP model to a discrete-time Markov decision
process (DTMDP) model based on the studies of the undiscounted
occupancy and occupation measures. We allow that the controlled
process is not necessarily absorbing, and the transition rates are
not necessarily separated from zero, and can be arbitrarily
unbounded; these features count for the main technical
difficulties in studying undiscounted CTMDP models.



\section{Introduction}
The present paper considers the constrained optimal control
problem with total undiscounted criteria for a continuous-time
Markov decision process (CTMDP) in Borel state and action spaces.

The majority of the previous literature on CTMDPs with the total
cost criteria focuses on the discounted model with a positive
constant discount factor; see e.g.,
\cite{Feinberg:2004,Feinberg:2012OHL,Guo:2009,GuoABP:2010,Kitaev:1995,Serfozo:1979,Piunovskiy:1998,ABPZY:20102,PiunovskiyZY:2012,PrietoOHL:2012book,Yushkevich:1980}.
In
\cite{Guo:2007,Guo:2009,GuoABP:2010,Piunovskiy:1998,ABPZY:20102,PrietoOHL:2012book},
 the convex analytic approach for constrained problems is
developed, whereas the dynamic programming approach for
unconstrained problems is studied in
\cite{Guo:2007,Guo:2009,PrietoOHL:2012book}.  The investigations
in
\cite{Guo:2007,Guo:2009,GuoABP:2010,Piunovskiy:1998,ABPZY:20102,PrietoOHL:2012book}
are based on the direct investigation of the continuous-time
models by using the Kolmogorov forward equations; for this, the
authors had to impose extra conditions bounding the growth of the
transition rates in the form of the existence of Lyapunov
functions.

Another method of investigation is based on the study of the
relation of the CTMDP problem and a DTMDP (discrete-time Markov
decision process) problem. Once the CTMDP problem is reduced to an
equivalent DTMDP problem, one can directly make use of the toolbox
of the better developed theory of DTMDPs
\cite{Altman:1999,BauerleRieder:2011,Bertsekas:1978,Piunovskiy:1997}
for the CTMDPs. This idea dates back to at least to the 1970s; see
Lippman \cite{Lippman:1975}, where the author applied the
uniformization technique to the reducing the CTMDP problem to a
DTMDP problem; see also \cite{Serfozo:1979}. However, the authors
of \cite{Lippman:1975,Serfozo:1979}, not only required the
transition rates to be uniformly bounded, also had to be
restricted to the class of deterministic stationary policies,
i.e., those that do not change actions between two consecutive
state transitions. These are also the standard setup for textbook
treatment of CTMDPs; see Chapter 11 of \cite{Puterman:1994}. The
situation becomes more complicated if one is allowed, as in the
present paper, to consider nonstationary policies, i.e., those
allowing the change in actions between two state transitions. In
this direction, Yushkevich \cite{Yushkevich:1980} firstly reduced
a discounted CTMDP model with nonstationary policies to a DTMDP
model. However, the action space of the induced DTMDP model is
more complicated; it is the space of measurable mappings, so that
in general a stationary policy in the DTMDP model corresponds to a
nonstationry policy for the original CTMDP model. A further
reduction of the induced DTMDP model to one with the same action
space as the original CTMDP model is possible after the
investigations of the dynamic programming (or say optimality)
equation for unconstrained problems; see Remark
\ref{GuoZhangAppendixRem} for greater details. Only unconstrained
problems were considered in \cite{Yushkevich:1980}, which also
assumed the transition rates in the CTMDP model to be uniformly
bounded.

In general, the reduction method based on the comparison of the
dynamic programming equation is more suitable for unconstrained
problem; see also \cite{PiunovskiyZY:2012}. Especially convenient
for dealing with constrained discounted CTMDP problems, Feinberg
\cite{Feinberg:2004,Feinberg:2012OHL} proposed a novel method of
reducing the CTMDP model to an equivalent DTMDP model in the same
action space based on the studies of the discounted occupancy
measures. (In fact, there is inconsistency in the use of
terminologies in \cite{Feinberg:2004,Feinberg:2012OHL}; the
occupation measure in \cite{Feinberg:2004} actually means the
occupancy measure in \cite{Feinberg:2012OHL} as well as the
present paper.) The original article \cite{Feinberg:2004} assumed
the transition rates of the CTMDP to be bounded; this condition is
completely withdrawn in the more recent extension
\cite{Feinberg:2012OHL}. Feinberg's reduction is valid without any
conditions so long the discount factor is positive.

All the aforementioned works are for discounted CTMDP models. The
present paper considers the total undiscounted CTMDP problem with
constraints. To the best of our knowledge, the theory for this
class of optimal control problems is currently underdeveloped,
despite that they would naturally find applications to e.g.,
epidemiology, where one aims at minimizing the total endemic time,
which does not have an obvious monetary interpretation for
discounting. There seems to be limited literature on this topic.
For unconstrained total undiscounted problem, Forwick et al
\cite{Forwick:2004} developed the dynamic programming approach,
and established the optimality equation, essentially following the
Yushevich's reduction method. For the constrained problem, the
authors of \cite{GuoMZhang:2013} developed the convex analytic
approach by studying directly the continuous-time model, but only
after imposing the extra conditions on the growth of the
transition rate and some strongly absorbing structure.

The objective of this paper is to study the constrained total
undiscounted CTMDP problem without the absorbing condition or any
condition on the growth of the transition rate. For DTMDPs, such
problems were acknowledged to be challenging in the survey
\cite{Borkar:2002} and were tackled only recently in
\cite{Dufour:2012}. Our original plan is to apply the Feinberg's
reduction method to the undiscounted case; we remark that the
Feinberg's reduction method is always applicable to discounted
CTMDP models without additional conditions. However, we notice
that the situation when the discount factor for the CTMDP model is
zero becomes significantly different and more delicate; indeed,
Example \ref{ZyChapterExample} below illustrates that without
additional conditions (in fact when the transition rate is not
separated from zero), it can happen that the performance vector of
the CTMDP problem under a nonstationary policy might not be
replicated by any performance vectors of the induced DTMDP
problem. It is thus natural to ask under what conditions does the
reduction method apply to the undiscounted CTMDP model. It is also
realized the studies of the occupancy measures alone are not
useful in general for the total undiscounted CTMDP models. (In
Section \ref{GuoZhangSec3} below we give a more detailed
discussion on these.) Different from the discounted case, we now
also need study the occupation measures, which are on the one
hand, more delicate because they are infinitely valued, and on the
other hand, are more suitable and convenient for constrained
problems.

Having said the above, the main contributions of the present paper
are as follows.
\begin{itemize}
\item[(a)] We provide the natural condition for the validity of
reducing the total undiscounted CTMDP model with constraints to a
DTMDP model. Our conditions are of the standard continuity and
compactness type, and allow the transition rates not necessarily
separated from zero on the one hand, and arbitrarily unbounded on
the other hand. No absorbing structure is assumed. The approach in
\cite{GuoMZhang:2013} are not applicable in this general setup.
Also note that the arguments in Feinberg
\cite{Feinberg:2004,Feinberg:2012OHL} are essentially based on the
presence of the positive discount factor; see Section
\ref{GuoZhangSec3} for greater details.

\item[(b)] We show the existence of an optimal stationary policy
out of the class of general (nonstationary) ones. It is arguable
that the solvability, as we confine ourselves to in this paper, is
an issue of core importance to be addressed first for any optimal
control problem.

\item [(c)] The paper is not a simple extension of the
uniformization technique for CTMDPs, as explained in the above.
Rather, our investigations are based on the studies of
undiscounted occupancy measures and occupation measures of the
CTMDP model, for which we incidentally obtain some properties of
independent interest.
\end{itemize}

The rest of this paper is organized as follows. We describe the
controlled process and state the concerned optimal control
problems in Section \ref{GuoZhangSec2}. In Sections
\ref{GuoZhangSec4} and \ref{GuoZhangSec5} we obtain some
properties of the occupancy and occupation measures, respectively.
In Section \ref{GuoZhangSec6} we establish the optimality results.
We end this paper with a conclusion in Section \ref{GuoZhangSec7}.
Some auxiliary statements and materials are presented in the
appendix.

\section{Optimal control problem statement}\label{GuoZhangSec2}
The objective of this section is to describe briefly the
controlled process similarly to
\cite{Kitaev:1986,Kitaev:1995,Piunovskiy:1998}, and the associated
optimal control problem of interest in this paper.
\bigskip

\par\noindent\textbf{Notations and conventions.} In what follows, $I$ stands for the indicator function, $\delta_{x}(\cdot)$
is the Dirac measure concentrated at $x,$ and ${\cal{B}}(X)$ is
the Borel $\sigma$-algebra of the topological space $X.$ The
abbreviation s.t. (resp., a.s.) stands for ``subject to'' (resp.,
``almost surely''). Below, unless stated otherwise, the term of
measurability is always understood in the Borel sense. Throughout
this article, we adopt the conventions of $\frac{0}{0}:=0,$
$0\cdot\infty:=0$ and $\frac{1}{0}:=+\infty.$

\subsection{Description of the CTMDP}

The primitives of a CTMDP model are the following elements $\{S,
A, q, \gamma\},$ where $S$ is a nonempty Borel state space, $A$ is
a nonempty Borel action space, $\gamma$ is a probability measure
on ${\cal B}(S)$ and represents the initial distribution, and $q$
stands for a signed kernel $q(dy|x,a)$ on ${\cal{B}}(S)$ given
$(x,a)\in S\times A$ such that
$\tilde{q}(\Gamma_S|x,a):=q(\Gamma_S\setminus\{x\}|x,a)\ge 0$ for
all $\Gamma_S\in{\cal{B}}(S).$ Throughout this article we assume
that $q(\cdot|x,a)$ is conservative and stable, i.e., $q(S|x,a)=0$
and $\bar{q}_x=\sup_{a\in A(x)}q_x(a)<\infty,$ where $q_x(a):=
-q(\{x\}|x,a).$ The signed kernel $q$ is often called the
transition rate. Throughout this article, $\bar{q}_x$ is allowed
to be arbitrarily unbounded in $x\in S$, unlike in
\cite{Guo:2007,GuoMZhang:2013,Piunovskiy:1998,ABPZY:20102}. In
line with \cite{Dufour:2012,Forwick:2004} and to fix ideas, we do
not consider the case of different admissible action spaces at
different states.

Let us take the sample space $\Omega$ by adjoining to the
countable product space $S\times((0,\infty)\times S)^\infty$ the
sequences of the form
$(x_0,\theta_1,\dots,\theta_n,x_n,\infty,x_\infty,\infty,x_\infty,\dots),$
where $x_0,x_1,\dots,x_n$ belongs to $S$,
$\theta_1,\dots,\theta_n$ belongs to $(0,\infty),$ and
$x_{\infty}\notin S$ is the isolated point. We equip $\Omega$ with
its Borel $\sigma$-algebra $\cal F$.

Let $t_0(\omega):=0=:\theta_0,$ and for each $n\geq 0$, and each
element $\omega:=(x_0,\theta_1,x_1,\theta_2,\dots)\in \Omega$, let
\begin{eqnarray*}
t_n(\omega)&:=&t_{n-1}(\omega)+\theta_n,
\end{eqnarray*}
and the limit point of the sequence $\{t_n\}$ is denoted by $
t_\infty(\omega):=\lim_{n\rightarrow\infty}t_n(\omega). $
Obviously, $t_n(\omega)$ are measurable mappings on $(\Omega,{\cal
F})$. In what follows, we often omit the argument $\omega\in
\Omega$ from the presentation for simplicity. Also, we regard
$x_n$ and $\theta_{n+1}$ as the coordinate variables, and note
that the pairs $\{t_n,x_n\}$ form a marked point process with the
internal history $\{{\cal F}_t\}_{t\ge 0},$ i.e., the filtration
generated by $\{t_n,x_n\}$; see Chapter 4 of \cite{Kitaev:1995}
for greater details. The marked point process $\{t_n,x_n\}$
defines the stochastic process on $(\Omega,{\cal F})$ of interest
$\{\xi_t,t\ge 0\}$ by
\begin{eqnarray}\label{GZCTMDPdefxit}
\xi_t=\sum_{n\ge 0}I\{t_n\le t<t_{n+1}\}x_n+I\{t_\infty\le
t\}x_\infty;
\end{eqnarray}
recall that $x_\infty$ is the isolated point. Below we denote
$S_{\infty}:=S\bigcup\{x_\infty\}.$

\begin{definition}
A (history-dependent) policy $\pi$ for the CTMDP is given by a
sequence $(\pi_n)$ such that, for each $n=1,2,\dots,$
$\pi_n(da|x_0,\theta_1,\dots,x_{n-1},s)$ is a stochastic kernel on
$A$, and for each $\omega=(x_0,\theta_1,x_1,\theta_2,\dots)\in
\Omega$, $t> 0,$
\begin{eqnarray*}
\pi(da|\omega,t)&=&I\{t\ge t_\infty\}\delta_{a_\infty}(da)+
\sum_{n=0}^\infty I\{t_n< t\le
t_{n+1}\}\pi_{n+1}(da|x_0,\theta_1,\dots,\theta_n,x_n, t-t_n),
\end{eqnarray*}
where $a_\infty\notin A$ is some isolated point. A policy
$\pi=(\pi_n)$  is called Markov if, with slight abuse of
notations, each of the stochastic kernels $\pi_n$ reads $
\pi_n(da|x_0,\theta_1,\dots,x_{n-1},s)=\pi_n(da|x_{n-1},s).$ A
Markov policy is further called deterministic if the stochastic
kernels $\pi_n(da|x_{n-1},s)$ all degenerate. A policy
$\pi=(\pi_n)$  is called stationary if, with slight abuse of
notations, each of the stochastic kernels $\pi_n$ reads $
\pi_n(da|x_0,\theta_1,\dots,x_{n-1},s)=\pi(da|x_{n-1}). $ A
stationary policy is further called deterministic if $
\pi_n(da|x_0,\theta_1,\dots,x_{n-1},s)=\delta_{f(x_{n-1})}(da) $
for some measurable mapping $f$ from $S$ to $A$.
\end{definition}
The class of all policies for the CTMDP model is denoted by $\Pi,$
and the class of all deterministic Markov policies for the CTMDP
model is denoted by $\Pi_{DM}.$

Under a policy $\pi:=(\pi_n)\in \Pi_{CTMDP}$, we define the
following random measure on $S\times (0,\infty)$
\begin{eqnarray*}
\nu^\pi(dt, dy)&:=& \int_A \tilde{q}(dy|\xi_{t-}(\omega),a)\pi(da|\omega,t)dt\nonumber\\
&=&\sum_{n\ge 0}\int_A
\tilde{q}(dy|x_n,a)\pi_{n+1}(da|x_0,\theta_1,\dots, \theta_n,
x_n,t-t_n)I\{t_n< t\le t_{n+1}\}dt
\end{eqnarray*}
with $q_{x_\infty}(a_\infty)=q(dy|x_\infty,a_\infty):=0. $ Then
there exists a unique probability measure ${P}^\pi_\gamma$ such
that
\begin{eqnarray*}
{P}_{\gamma}^\pi(x_0\in dx)=\gamma(dx),
\end{eqnarray*}
and with
respect to $P_\gamma^\pi,$ $\nu^\pi$ is the dual predictable
projection of the random measure associated with the marked point
process $\{t_n,x_n\}$; see \cite{Jacod:1975,Kitaev:1995}. The
process $\{\xi_t\}$ defined by (\ref{GZCTMDPdefxit}) under the
probability measure ${}{P}_\gamma^\pi$ is called a CTMDP. Below,
when $\gamma(\cdot)$ is a Dirac measure concentrated at $x\in S,$
we use the denotation ${}{P}_x^\pi.$ Expectations with respect to
${}{P}_\gamma^\pi$ and ${}{P}_x^\pi$ are denoted as
${}{E}_{\gamma}^\pi$ and ${}{E}_{x}^\pi,$ respectively.

Under the probability measure $P_\gamma^\pi,$ the system dynamics
of the CTMDP can be described as follows. The initial state $x_0$
has the distribution given by $\gamma$, the sojourn time
$\theta_{n+1}$ has the (conditional) tail function given by
\begin{eqnarray*}
P_\gamma^\pi(\theta_{n+1}\ge t
|x_0,\theta_1,\dots,x_n)=e^{-\int_0^t
\int_{A}q_{x_n}(a)\pi_{n+1}(da|x_0,\theta_1,\dots,x_n,s)ds},
\end{eqnarray*}
and upon a jump, the (conditional) distribution of the next state
$x_{n+1}$ is given by
\begin{eqnarray*}
P_\gamma^\pi(x_{n+1}\in
\Gamma|x_0,\theta_1,\dots,x_n,\theta_{n+1})=\frac{\int_A
\tilde{q}(\Gamma|x_n,a)\pi_{n+1}(da|x_0,\theta_1,\dots,x_n,\theta_{n+1})}{\int_A
q_{x_n}(a)\pi_{n+1}(da| x_0,\theta_1,\dots,x_n,\theta_{n+1})}
\end{eqnarray*}
for each $\Gamma\in {\cal B}(S).$ Here and below we formally put
$\pi_{n+1}(\{a_\infty\}|x_0,\theta_1,\dots,x_n,\infty):=1$ with
$a_\infty\notin A$ being the isolated point, so that $\int_A
\tilde{q}(\Gamma|x_n,a)\pi_{n+1}(da|x_0,\theta_1,\dots,x_n,\infty):=0$
for each $\Gamma\in {\cal B}(S)$. Also recall the convention of
$\frac{0}{0}:=0,$ so that
\begin{eqnarray*}
P_\gamma^\pi(x_{n+1}=x_\infty|x_0,\theta_1,\dots,x_n,\theta_{n+1})=1-P_\gamma^\pi(x_{n+1}\in
S|x_0,\theta_1,\dots,x_n, \theta_{n+1}).
\end{eqnarray*}

In what follows, when it is not necessary to emphasize the initial
distribution $\gamma$, we also say that $\{S,A,q\}$ is our CTMDP
model.

\subsection{Description of the concerned optimal control problem}

Let $N\in\{1,2,\dots\}$ be fixed. Consider the nonnegative
measurable functions $c_i(x,a)\ge 0$ with $i=0,1,\dots,N$ from
$S\times A$ to $[0,\infty)$ as the cost rates.  We formally put
$c_i(x_\infty,a):=0$ for each $i=0,1,2,\dots,N.$

In this paper, we study the following optimal control problem:
\begin{eqnarray}\label{ZyChapter1}
&&E_\gamma^\pi\left[\int_0^\infty \int_A
c_0(\xi_t,a)\pi(da|\omega,t)dt\right]\rightarrow \min_{\pi\in
\Pi}\nonumber \\
&s.t.& E_\gamma^\pi\left[\int_0^\infty \int_A
c_j(\xi_t,a)\pi(da|\omega,t)dt\right]\le
d_j,~\forall~j=1,2,\dots,N,
\end{eqnarray}
where for each $j=1,2,\dots,N,$ $d_j\in[0,\infty)$ is the fixed
constraint constant.

A policy $\pi\in \Pi$ is called feasible for problem
(\ref{ZyChapter1}) if
\begin{eqnarray*}
E_\gamma^\pi\left[\int_0^\infty \int_A
c_j(\xi_t,a)\pi(da|\omega,t)dt \right]\le d_j,~j=1,2,\dots,N.
\end{eqnarray*}
Let $\Pi_F$ be the class of feasible policies. Then the value of
problem (\ref{ZyChapter1}) is denoted as
\begin{eqnarray*}
V_c(\gamma):=\inf_{\pi\in \Pi_F} E_\gamma^\pi\left[\int_0^\infty
\int_A c_0(\xi_t,a)\pi(da|\omega,t)dt\right].
\end{eqnarray*}

A feasible policy $\pi$ for problem (\ref{ZyChapter1}) is called
to be with a finite value if
\begin{eqnarray*}
E_\gamma^\pi\left[\int_0^\infty \int_A
c_0(\xi_t,a)\pi(da|\omega,t)dt \right]<\infty.
\end{eqnarray*}

Finally, a policy $\pi^\ast\in \Pi$ is called optimal for the
(constrained) CTMDP problem (\ref{ZyChapter1}) if it holds that
\begin{eqnarray*}
\inf_{\pi\in \Pi_F} E_\gamma^\pi\left[\int_0^\infty \int_A
c_0(\xi_t,a)\pi(da|\omega,t)dt\right]=
E_\gamma^{\pi^\ast}\left[\int_0^\infty \int_A
c_0(\xi_t,a)\pi^\ast(da|\omega,t)dt\right].
\end{eqnarray*}

\section{Facts about the discounted CTMDP problem and discussions}
\label{GuoZhangSec3}

The purpose of this section is to (a) present some relevant
results about the $\alpha$-discounted problem for the CTMDP model
$\{S,A,q,\gamma\}$, which are used in the subsequent
investigations for our undiscounted CTMDP problem
(\ref{ZyChapter1}); and (b) demonstrate the significant difference
between the discounted and the undiscounted CTMDP problems, and
illustrate that the undiscounted problem is more delicate, which
thus clarifies the contribution of the present paper; see Example
\ref{ZyChapterExample} and the discussion following it.

In his well written articles
\cite{Feinberg:2004,Feinberg:2012OHL}, Professor Feinberg
considered the following constrained discounted optimal control
problem for the CTMDP model $\{S,A,q,\gamma\}$
\begin{eqnarray}\label{FeinbergChapter1}
&&E_\gamma^\pi\left[\int_0^\infty e^{-\alpha t}\int_A
c_0(\xi_t,a)\pi(da|\omega,t)dt\right]\rightarrow \min_{\pi\in
\Pi}\nonumber \\
&s.t.& E_\gamma^\pi\left[\int_0^\infty  e^{-\alpha t} \int_A
c_j(\xi_t,a)\pi(da|\omega,t)dt\right]\le
d_j,~\forall~j=1,2,\dots,N,
\end{eqnarray}
where $d_j\in \mathbb{R}$ for each $j=1,2,\dots,N,$ where the
finite constant $\alpha>0$ is a fixed discount factor. The
investigations in \cite{Feinberg:2004,Feinberg:2012OHL} are based
on the study of the so-called $\alpha$-discounted occupancy
measures, firstly introduced therein, which we recall as follows.

\begin{definition}\label{FeinbergChapterDef1}
For each $n=0,1,\dots,$ and (finite) constant $\alpha>0$, the
$\alpha$-discounted occupancy measure of the policy $\pi\in \Pi$
for the CTMDP model $\{S,A,q,\gamma\}$ is a measure
$M_{\gamma,\alpha}^{n,\pi}$ on ${\cal B}(S\times A)$ defined by
for each $\Gamma_S \in {\cal B}(S)$ and $\Gamma_A\in {\cal B}(A)$,
\begin{eqnarray*}
M_{\gamma,\alpha}^{n,\pi}(\Gamma_S\times \Gamma_A):=
E_\gamma^\pi\left[\int_{t_n}^{t_{n+1}} e^{-\alpha t} I\{\xi_t\in
\Gamma_S\} \int_{\Gamma_A} (\alpha+q_{\xi_t}(a))\pi(da|\omega,t)dt
\right].
\end{eqnarray*}
\end{definition}

Professor Feinberg noticed that there is a close relationship
between the ($\alpha$-discounted) occupancy measure for the CTMDP
model $\{S,A,q,\gamma\}$ and the marginal distribution of
$(X_n,A_{n+1})$ of the DTMDP model
$\{S_\infty,A,p_\alpha,\gamma\}$, where the transition probability
$p_\alpha$ is defined by for each $\Gamma_S\in {\cal B}(S)$,
\begin{eqnarray*}
p_\alpha(\Gamma_S|x,a)=\frac{\tilde{q}(\Gamma_S|x,a)}{\alpha+q_x(a)},~\forall~x\in
S,~a\in A
\end{eqnarray*}
and
\begin{eqnarray*}
p_\alpha(\Gamma_S|x_\infty,a)=0,~\forall~a\in A.
\end{eqnarray*}
Recall that $S_\infty=S\bigcup\{x_\infty\}$ with $x_\infty\notin
S$ being the isolated point. Under each policy $\sigma$ for the
DTMDP model $\{S_\infty, A,p_\alpha,\gamma\}$, let the
corresponding strategic measure be denoted by
$\textbf{P}_\gamma^{\alpha,\sigma}$. The expectation taken with
respect to $\textbf{P}_\gamma^{\alpha,\sigma}$ is written as
$\textbf{E}_\gamma^{\alpha,\sigma}.$

The next statement is established in \cite{Feinberg:2012OHL}.
\begin{proposition}\label{FeinbergChapterTheorem1}
The following assertions hold for each $\Gamma_S\in {\cal B}(S)$
and $\Gamma_A\in {\cal B}(A)$.

\par\noindent(a) For each
policy $\pi\in\Pi$ for the CTMDP model, there is a (possibly)
randomized Markov policy $\sigma^M$ for the DTMDP model
$\{S_\infty, A,p_\alpha\}$ such that
\begin{eqnarray*}\label{FeinbergChapter6}
M_{\gamma,\alpha}^{n,\pi}(\Gamma_S\times
\Gamma_A)=\textbf{P}_\gamma^{\alpha,\sigma^M}(X_{n}\in
\Gamma_S,~A_{n+1}\in \Gamma_A),~\forall~n=0,1,\dots.
\end{eqnarray*}
\par\noindent(b) For each randomized Markov policy $\sigma^M$ for the DTMDP model $\{S_\infty,A,p_\alpha\}$, there exists
a Markov policy $\pi^M$ for the CTMDP model such that
\begin{eqnarray*}
M_{\gamma,\alpha}^{n,\pi^M}(\Gamma_S\times
\Gamma_A)=\textbf{P}_\gamma^{\alpha,\sigma^M}(X_{n}\in
\Gamma_S,~A_{n+1}\in \Gamma_A),~\forall~n=0,1,\dots.
\end{eqnarray*}
\end{proposition}

Proposition \ref{FeinbergChapterTheorem1} shows that for each
$\pi\in \Pi$ for the CTMDP model $\{S,A,q,\gamma\},$ there exists
some policy $\sigma$ for the DTMDP model
$\{S_\infty,A,p_\alpha,\gamma\}$ such that
\begin{eqnarray}\label{ZyChapter73}
 \sum_{n=0}^\infty M_{\gamma,\alpha}^{n,\pi}(dx\times da)=\sum_{n=0}^\infty \textbf{P}_\gamma^{\alpha,\sigma}(X_n\in
dx,~A_{n+1}\in da).
\end{eqnarray}
On the opposite direction, for each Markov policy $\sigma^M$ for
the DTMDP model $\{S_\infty,A,p_\alpha,\gamma\}$, there exists
some policy $\pi$ for the CTMDP model $\{S,A,q,\gamma\},$ such
that
\begin{eqnarray*}
 \sum_{n=0}^\infty M_{\gamma,\alpha}^{n,\pi}(dx\times da)=\sum_{n=0}^\infty \textbf{P}_\gamma^{\alpha,\sigma^M}(X_n\in
dx,~A_{n+1}\in da).
\end{eqnarray*}
Furthermore, for each $\sigma$ for the DTMDP model
$\{S_\infty,A,p_\alpha,\gamma\}$, there is a (possibly) randomized
Markov policy $\sigma^M$ for the DTMDP model such that
\begin{eqnarray*}
\textbf{E}_\gamma^{\alpha,\sigma}\left[\sum_{n=0}^\infty
\frac{c_i(X_n,A_{n+1})}{\alpha+q_{X_n}(A_{n+1})}
 \right]=\textbf{E}_\gamma^{\alpha,\sigma^M}\left[\sum_{n=0}^\infty \frac{c_i(X_n,A_{n+1})}{\alpha+q_{X_n}(A_{n+1})}
 \right];
 \end{eqnarray*}
this is due to the well known Derman-Strauch lemma
\cite{Derman:1966}; see also Lemma 2 of Piunovskiy
\cite{Piunovskiy:1997}. Consequently, Proposition
\ref{FeinbergChapterTheorem1} shows that the $\alpha$-discounted
CTMDP problem (\ref{FeinbergChapter1}) can be reduced to the
following DTMDP problem for the model
$\{S_\infty,A,p_\alpha,\gamma\}$
\begin{eqnarray*}
&&\textbf{E}_\gamma^{\alpha,\sigma}\left[\sum_{n=0}^\infty
\frac{c_0(X_n,A_{n+1})}{\alpha+q_{X_n}(A_{n+1})}\right]\rightarrow\min_{\sigma}
\\
&s.t.& \textbf{E}_\gamma^{\alpha,\sigma}\left[\sum_{n=0}^\infty
\frac{c_j(X_n,A_{n+1})}{\alpha+q_{X_n}(A_{n+1})}\right]\le
d_j,~j=1,2,\dots,N.
\end{eqnarray*}
(Recall that $c_i(x_\infty,a):=0$ for each $a\in A$.)

Here and below by reduction is meant that both problems have the
same value, and if an optimal policy exists for one problem, so
does an optimal policy for the other problem.

We emphasize that this reduction for the $\alpha$-discounted CTMDP
problem is possible without any extra conditions being imposed on
the CTMDP model, so long $\alpha>0.$

It is natural to ask whether the reduction is possible for the
case of $\alpha=0;$ i.e., whether the CTMDP problem
(\ref{ZyChapter1}) can be reduced to the following problem
\begin{eqnarray}\label{ZyChapter63}
&&\textbf{E}_\gamma^\sigma\left[\sum_{n=0}^\infty
\frac{c_0(X_n,A_{n+1})}{q_{X_n}(A_{n+1})}\right]\rightarrow
\min_{\sigma}\nonumber\\
&s.t.&\textbf{E}_\gamma^\sigma\left[\sum_{n=0}^\infty
\frac{c_j(X_n,A_{n+1})}{q_{X_n}(A_{n+1})}\right]\le
d_j,~\forall~j=1,2,\dots,N.
\end{eqnarray}
for the DTMDP model $\{S_\infty,A,p,\gamma\},$ where the
transition probability $p$ being defined by for each $\Gamma_S\in
{\cal B}(S)$,
\begin{eqnarray}\label{ZyChapter64}
p(\Gamma_S|x,a)=\frac{\tilde{q}(\Gamma_S|x,a)}{q_x(a)},~\forall~x\in
S,~a\in A
\end{eqnarray}
and
\begin{eqnarray}\label{ZyChapter65}
p(\Gamma_S|x_\infty,a)=0.
\end{eqnarray}
(Recall that $\frac{0}{0}:=0.$) As before, the controlled and
controlling processes for the DTMDP model
$\{S_\infty,A,p,\gamma\}$ are denoted by $\{X_n\}$ and $\{A_n\}$;
$\textbf{P}_\gamma^\pi$ denotes the strategic measure under the
policy $\sigma$ for this DTMDP model with the corresponding
expectation $\textbf{E}_\gamma^\sigma$.

We remark that since $c_i(x_\infty,a)=0$ and
$p(dy|x_\infty,a)=\delta_{x_\infty}(dy)$ for each $a\in A,$ the
definition of a policy $\sigma$ at the current state $x_\infty$
for the DTMDP model $\{S_\infty,A,p,\gamma\}$ is not important for
its performance as far as problem (\ref{ZyChapter63}) is
concerned, and so we do not specify it in what follows.

The next example shows that the answer to the above question is
negative in general.

\begin{example}\label{ZyChapterExample}
Consider the CTMDP model with $S=\{0,1\}$, $A=[0,\infty)$,
$q_1(a)=q(\{2\}|1,a)=e^{-a}$, $q_2(a)=0$ for each $a\in A,$ and
$\gamma(\{1\})=1.$ Let $N=1,$ and $c_0(1,a)=e^{-a}$, $c_0(2,a)=0$
for each $a\in A,$ and $c_1(x,a)=0$ for each $x\in S$ and $a\in
A.$ Let $d_1>0,$ so that any policy is feasible for the CTMDP
problem (\ref{ZyChapter1}). Let us fix a policy $\pi$ defined by
for each $a\in A,$
\begin{eqnarray*}
\pi(\{a\}|\omega,t)=\pi_{0}(\{a\}|x,t)=I\{a=t\},
\end{eqnarray*}
so that
\begin{eqnarray*}
\int_{A}q_1(a)\pi(da|1,t)=\int_A c_0(1,a)\pi(da|1,t)=e^{-t}.
\end{eqnarray*}
Then under this policy $\pi$, we see
\begin{eqnarray*}
E_\gamma^\pi\left[\int_0^\infty \int_A
c_0(\xi_t,a)\pi(da|\omega,t)dt\right]= E_\gamma^\pi\left[
\int_0^{\theta_1}e^{-t} dt\right]<\int_0^\infty e^{-t}dt=1,
\end{eqnarray*}
where the third equality is due to the fact
$P_\gamma^\pi(\theta_1=\infty)=e^{-1}<1$. On the other hand, since
$\frac{c_0(1,a)}{q_1(a)}=1$ for each $a\in A$, we have that under
each policy $\sigma$ for the DTMDP model $\{S_\infty,A,p,\gamma\}$
\begin{eqnarray*}
\textbf{E}_\gamma^\sigma\left[\sum_{n=0}^\infty
\frac{c_0(X_n,A_{n+1})}{q_{X_n}(A_{n+1})}\right]= 1.
\end{eqnarray*}
In summary, each policy for the DTMDP $\{S_\infty,A,p,\gamma\}$
model would be optimal for problem (\ref{ZyChapter63}) with the
(optimal) value being $1$, whereas the value for the CTMDP problem
(\ref{ZyChapter1}) is strictly smaller than $1$. Hence, the CTMDP
problem (\ref{ZyChapter1}) cannot be reduced to the DTMDP problem
(\ref{ZyChapter63}).
\end{example}

It is also clear that Proposition \ref{FeinbergChapterTheorem1}
does not hold in general when $\alpha=0;$ see also Remark
\ref{ZyChapterRemark1} below.

An objective of the present paper is to provide weak and natural
conditions under which the reduction of the CTMDP problem
(\ref{ZyChapter1}) to the DTMDP problem (\ref{ZyChapter63}) is
possible. To this end, apart from studying the (undiscounted)
occupancy measures (see Definition \ref{ZyChapterDefinition1}), we
also need investigate the (undiscounted) occupation measures (see
Definition \ref{ForwickChapterDef1}) for the CTMDP model
$\{S,A,q,\gamma\}$, for which some properties are to be obtained.
The occupation measure is more delicate for studies because it is
infinitely valued, whereas the occupancy measure is always finite;
see (\ref{ZyChapter11}) below. Finally, under our conditions, we
obtain the existence of an optimal stationary policy for the CTMDP
problem (\ref{ZyChapter1}). It is arguable that the solvability,
as we confine ourselves to in this paper, is an issue of core
importance to be addressed for any optimal control problem.

\section{Occupancy measure}\label{GuoZhangSec4}

The objective in this section is to obtain a partial version of
Theorem \ref{FeinbergChapterTheorem1}(a); see Theorem
\ref{ZyChapterTheorem2} below. This statement is needed in the
subsequent sections.

\begin{definition}\label{ZyChapterDefinition1}
For each $n=0,1,\dots,$ the (undiscounted) occupancy measure of
the policy $\pi\in \Pi$ for the CTMDP model $\{S,A,q,\gamma\}$ is
a measure $M_{\gamma}^{n,\pi}$ on ${\cal B}(S\times A)$ defined by
for each $\Gamma_S \in {\cal B}(S)$ and $\Gamma_A\in {\cal B}(A)$,
\begin{eqnarray}\label{ZyChapter66}
M_{\gamma}^{n,\pi}(\Gamma_S\times \Gamma_A):=
E_\gamma^\pi\left[\int_{t_n}^{t_{n+1}} I\{\xi_t\in \Gamma_S\}
\int_{\Gamma_A} q_{\xi_t}(a)\pi(da|\omega,t)dt \right].
\end{eqnarray}
\end{definition}

\begin{condition}\label{BookWC1}
\par\noindent (a) The space $A$ is compact.
\par\noindent(b) For each bounded continuous
function $f(x)$ on $S$,
$\int_{S}f(y)\frac{\tilde{q}(dy|x,a)}{q_x(a)}$ is continuous in
$(x,a)\in S\times A.$
\par\noindent(c) $q_x(a)$ is continuous in $(x,a)\in S\times A.$
\par\noindent(d) For each $i=0,1,\dots,N,$
$c_i(x,a)$ is lower semicontinuous in $(x,a)\in S\times A.$
\end{condition}

Let us introduce the following sets
\begin{eqnarray}\label{ZyChapter19}
S_1&:=&\left\{x\in S: \inf_{a\in A} q_x(a)=0,~ \inf_{a\in
A}\left(q_x(a)+\sum_{i=0}^N c_i(x,a)\right)>0\right\},\nonumber\\
\hat{S}_1&:=&\left\{x\in S_1: \sup_{a\in A}q_x(a)=0\right\},\nonumber\\
S_2&:=&\left\{x\in S: \inf_{a\in A} \left(q_x(a)+\sum_{i=0}^N c_i(x,a)\right)=0\right\},\nonumber\\
S_3&:=&\left\{x\in S: \inf_{a\in A} q_x(a)>0\right\}.
\end{eqnarray}
Under Condition \ref{BookWC1}, the above four sets are all
measurable, by Proposition 7.32 of \cite{Bertsekas:1978} and Lemma
\ref{ForwickChapterLem2}. Furthermore, $S_1,$ $S_2$ and $S_3$ are
disjoint and satisfy
\begin{eqnarray*}
S=S_1\bigcup S_2\bigcup S_3.
\end{eqnarray*}

Let us also denote for each $x\in S$,
\begin{eqnarray}\label{ZyChapter24}
B(x):=\{a\in A:~q_x(a)=0\},
\end{eqnarray}
which is compact under Condition \ref{BookWC1}.

\begin{lemma}\label{ZyChapterLem5}
Suppose Condition \ref{BookWC1} is satisfied. Consider a feasible
policy $\pi=(\pi_n)\in \Pi$ with a finite value for the CTMDP
problem (\ref{ZyChapter1}). Then for each $n=0,1,\dots$,
\begin{eqnarray*}
&&P_\gamma^\pi(x_n\in S_1\setminus \hat{S}_1)\\
&=&E_\gamma^\pi\left[I\{x_n\in S_1\setminus
\hat{S}_1\}\left.P_\gamma^\pi\left(\int_0^\infty \int_{A\setminus
B(x_{n})}q_{x_{n}}(a)\pi_{n+1}(da|x_0,\theta_1,\dots,x_n,s)
ds=\infty\right|x_n\right)\right].
\end{eqnarray*}
\end{lemma}

\par\noindent\textit{Proof.} It holds that for each $n=0,1,\dots,$
\begin{eqnarray}\label{ZyChapter10}
&&\infty>E_\gamma^\pi\left[\int_{t_n}^{t_{n+1}} \int
_{(A\setminus B(x_n))\bigcup B(x_n)}\sum_{i=0}^N c_i(x_n,a)\pi_{n+1}(da|x_0,\dots,x_n,t)dt\right]\nonumber\\
&\ge& E_\gamma^\pi\left[  \int_{t_n}^{t_{n+1}}I\{x_n\in
S_1\setminus\hat{S}_1\} \min _{a\in
B(x_n)}\left\{\sum_{i=0}^N c_i(x_n,a)\right\}ds\right]\nonumber \\
&=&E_\gamma^\pi\left[ I\{x_n\in S_1\setminus\hat{S}_1\} \min
_{a\in
  B(x_n)}\left\{\sum_{i=0}^N c_i(x_n,a)\right\} \int_0^\infty e^{-\int_0^t \int_{A}
  q_{x_n}(a)\pi_{n+1}(da|x_0,\dots,x_n,s)ds}dt\right].
\end{eqnarray}
If the statement of the lemma does not hold, then there is some
$n=0,1,\dots$ such that
\begin{eqnarray*}
P_\gamma^\pi\left( x_n\in S_1\setminus\hat{S}_1,~
 \int_0^\infty \int_{A}
  q_{x_n}(a)\pi_{n+1}(da|x_0,\theta_1,\dots,x_n,s)ds <\infty\right)>0,
\end{eqnarray*}
and thus
\begin{eqnarray*}
P_\gamma^\pi( x_n\in S_1\setminus\hat{S}_1,~\int_0^\infty
e^{-\int_0^t \int_{A}
  q_{x_n}(a)\pi_{n+1}(da|x_0,\theta_1,\dots,x_n,s)ds}dt=\infty)>0.
\end{eqnarray*}
This implies
\begin{eqnarray*}
E_\gamma^\pi\left[ I\{x_n\in S_1\setminus\hat{S}_1\} \min _{a\in
  B(x_n)}\left\{\sum_{i=0}^N c_i(x_n,a)\right\}\int_0^\infty e^{-\int_0^t \int_{A}
  q_{x_n}(a)\pi_{n+1}(da|x_0,\dots,x_n,s)ds}dt\right]
=\infty,
\end{eqnarray*}
where the last equality follows from the fact that $\min _{a\in
  B(x)}\left\{\sum_{i=0}^N c_i(x,a)\right\}>0$ for each $x\in S_1.$
This contradicts (\ref{ZyChapter10}). $\hfill\Box$
\bigskip

\begin{definition}
For each fixed $n=0,1,\dots,$ and policy $\pi=(\pi_n)\in \Pi,$ we
define a measure $m_{\gamma,n}^\pi(dx\times da)$ for the CTMDP
model $\{S,A,q,\gamma\}$ on ${\cal B}(S\times A)$ by
\begin{eqnarray}\label{ZyChapter17}
m_{\gamma,n}^\pi(\Gamma_S\times
\Gamma_A):=E_\gamma^\pi\left[\int_{t_n}^{t_{n+1}}  I\{x_n\in
\Gamma_S\}\pi_{n+1}(\Gamma_A|x_0,\theta_1,\cdots,x_n,t-t_n)dt
\right]
\end{eqnarray}
for each $\Gamma_S\in {\cal B}(S)$ and $\Gamma_A\in {\cal B}(A)$.
\end{definition}
Evidently, for each $n=0,1,\dots,$ and $\Gamma_S\in {\cal B}(S)$,
$m_{\gamma,n}^\pi(\Gamma_S\times A)>0$ if and only if
$P_\gamma^\pi(x_n\in \Gamma_S)>0.$

\begin{lemma}\label{ZyChapterLem6}
Suppose Condition \ref{BookWC1} is satisfied. Consider a feasible
policy $\pi\in \Pi$ with a finite value for the CTMDP problem
(\ref{ZyChapter1}). Then it holds that
\begin{eqnarray*}
\int_{S\times A}\tilde{q}(\hat{S}_1|x,a)m_{\gamma,n}^\pi(dx\times
da)=0,
\end{eqnarray*}
and
\begin{eqnarray}\label{ZyChapter58}
m_{\gamma,n}^\pi(\hat{S}_1\times A)=0
\end{eqnarray}
for each $n=0,1,\dots.$ In particular,
\begin{eqnarray*}
\gamma(\hat{S}_1)=0.
\end{eqnarray*}
\end{lemma}
\par\noindent\textit{Proof.}
Suppose for contradiction that $m_{\gamma,n}^\pi(\hat{S}_1\times
A)>0$ for some $n.$ Then similarly to the proof of Lemma
\ref{ZyChapterLem5}, one can establish the following
contradiction;
\begin{eqnarray*}
 \infty> E_\gamma^\pi\left[\int_{t_n}^{t_{n+1}} \int _A
 \sum_{i=1}^N c_i(x_n,a) \pi(da|\omega,t)dt \right] \ge
E_\gamma^\pi\left[\int_{t_n}^{t_{n+1}} \min_{a\in A}\left\{
\sum_{i=1}^N c_i(x_n,a)\right\}I\{x_n\in \hat{S}_1\} dt
\right]=\infty.
\end{eqnarray*}
As a result, $m_{\gamma,n}^\pi(\hat{S}_1\times A)=0$ for each
$n=0,1,\dots.$ In particular, $m_{\gamma,0}^\pi(\hat{S}_1\times
A)=0$, and thus $P_\gamma^\pi(x_0\in
\hat{S}_1)=\gamma(\hat{S}_1)=0.$ It remains to prove
$\int_{S\times A}\tilde{q}(\hat S_1|x,a) m_{\gamma,n}^\pi(dx\times
da)=0$ for each $n=0,1,\dots$. If this is not true,  then it
follows from the definition of $m^\pi_{\gamma,n}(dx\times da)$
that for some $n=0,1,\dots,$
\begin{eqnarray*}
0&<&E_\gamma^\pi\left[\int_0^{\theta_{n+1}} \int_A
\tilde{q}(\hat{S}_1|x_n,a)\pi_{n+1}(da|x_0,\theta_1,\dots,x_n,t)
dt \right]\\
&=&E_\gamma^\pi\left[\int_0^\infty  \int_A
\tilde{q}(\hat{S}_1|x_n,a)\pi_{n+1}(da|x_0,\theta_1,\dots,x_n,t)
e^{-\int_0^t \int_A
q_{x_n}(a)\pi_{n+1}(da|x_0,\theta_1,\dots,x_n,s)ds }dt \right],
\end{eqnarray*}
which implies that $P_\gamma^\pi(x_{n+1}\in \hat{S}_1)>0$ by the
construction of the CTMDP; see (2) of \cite{Piunovskiy:1998}. This
leads to the contradiction against the fact that
$m_{\gamma,n+1}^\pi(\hat{S}_1\times A)>0$ as established earlier.
$\hfill\Box$

\begin{definition}
Let $f^\ast$ be a fixed measurable mapping from $S$ to $A$ such
that
\begin{eqnarray}\label{ZyChapter14}
0=\inf_{a\in A} \left\{\sum_{i=0}^N c_i(x,a)+q_x(a)\right\}=
\sum_{i=0}^N c_i(x,f^\ast(x))+q_x(f^\ast(x))
\end{eqnarray}
for each $x\in S_2$ whenever $S_2$ is nonempty.
\end{definition}
The existence of such a mapping is guaranteed by Proposition 7.33
of Bertsekas and Shreve \cite{Bertsekas:1978} under Condition
\ref{BookWC1}.

\begin{theorem}\label{ZyChapterTheorem2}
Suppose Condition \ref{BookWC1} is satisfied. Consider a feasible
policy $\pi=(\pi_n)\in \Pi$ with a finite value for the CTMDP
problem (\ref{ZyChapter1}) such that
\begin{eqnarray}\label{ZyChapter15}
\pi_{n+1}(da|x_0,\theta_1,\dots,x_n,s)=\delta_{f^\ast(x_n)}(da)
\end{eqnarray}
whenever $x_n\in S_2.$ Then there is a Markov policy $\sigma$ for
the DTMDP $\{S_\infty, A,p,\gamma\}$ such that for each
$n=0,1,\dots,$
\begin{eqnarray}\label{ZyChapter60}
\sigma_{n+1}(da|x)=\delta_{f^\ast(x)}(da)
\end{eqnarray}
for each $x\in S_2$ (if $S_2\ne \emptyset$), and
\begin{eqnarray}\label{ZyChapter13}
  M_{\gamma}^{n,\pi}(\Gamma_S\times \Gamma_A)=\textbf{P}_\gamma^\sigma(X_n\in \Gamma_S,~A_{n+1}\in
  \Gamma_A),~\forall~\Gamma_S\in {\cal B}(S\setminus S_2),~ \Gamma_A\in {\cal
  B}(A).
\end{eqnarray}
\end{theorem}

\par\noindent\textit{Proof.}
For each $\Gamma_S\in {\cal B}(S),$
\begin{eqnarray}\label{ZyChapter11}
&&\left.M_{\gamma}^{n,\pi}(\Gamma_S\times
A)=E_\gamma^\pi\left[E_\gamma^\pi
\left[\int_{t_n}^{t_{n+1}}I\{\xi_t\in \Gamma_S\} \int_{A}
q_{\xi_t}(a)\pi(da|\omega,t)dt\right|x_0,\theta_1,\dots,x_n
\right] \right]\nonumber\\
&=& \left.E_\gamma^\pi\left[I\{x_n\in \Gamma_S\}\right.
\left.E_\gamma^\pi\left[ \int_0^{\theta_{n+1}} \int_{A}
q_{x_n}(a)\pi_{n+1}(da|x_0,\theta_1,\dots,x_n,t)dt
\right|x_0,\theta_1,\dots,x_n \right]\right]\nonumber\\
&=& E_\gamma^\pi\left[I\{x_n\in \Gamma_S\} \int_0^\infty \int_{A}
q_{x_n}(a)\pi_{n+1}(da|x_0,\theta_1,\dots,x_n,t)e^{-\int_0^t
\int_A q_{x_n}(a)\pi_{n+1}(da|x_0,\theta_1,\dots,x_n,s)ds}dt
 \right]\nonumber\\
&=& E_\gamma^\pi\left[ I\{x_n\in \Gamma_S\}
\left(1-e^{-\int_0^\infty \int_A
q_{x_n}(a)\pi_{n+1}(da|x_0,\theta_1,\dots,x_n,s)ds}\right)\right]\le
1.
\end{eqnarray}
Then for each $n=0,1,\dots,$ one can refer to Corollary 7.27.2 of
\cite{Bertsekas:1978} or Proposition D.8 of
\cite{Hernandez-Lerma:1996} for the existence of a stochastic
kernel $\sigma_{n+1}(da|x)$ such that
\begin{eqnarray}\label{ZyChapter12}
M_{\gamma}^{n,\pi}(dx\times da)=M_{\gamma}^{n,\pi}(dx\times
A)\sigma_{n+1}(da|x).
\end{eqnarray}
on ${\cal B}(S\times A)$, and (\ref{ZyChapter60}) holds, where the
last assertion is true because $M_\gamma^{n,\pi}(S_2\times A)=0$
by (\ref{ZyChapter15}), (\ref{ZyChapter14}) and
(\ref{ZyChapter66}). Let $\sigma=(\sigma_n)$ be the Markov policy
for the DTMDP model $\{S_\infty,A,p,\gamma\}$ defined by this
sequence of stochastic kernels.

Consider the case of $n=0.$ Then for each $\Gamma_S\in {\cal
B}(S\setminus S_2),$
\begin{eqnarray}\label{ZyChapter38}
&&M_{\gamma}^{0,\pi}(\Gamma_S\times A)=E_\gamma^\pi\left[
I\{x_0\in \Gamma_S\} \left(1-e^{-\int_0^\infty \int_A
q_{x_0}(a)\pi_{1}(da|x_0,s)ds}\right)\right]\nonumber\\
&=&\gamma(\Gamma_S)=\textbf{P}_{\gamma}^\sigma(X_0\in \Gamma_S),
\end{eqnarray}
where the first equality is by (\ref{ZyChapter11}), the second
equality follows from Lemma \ref{ZyChapterLem6} in case
$\Gamma_S\subseteq \hat{S}_1$, from Lemma \ref{ZyChapterLem5} in
case $\Gamma_S\subseteq S_1\setminus \hat{S}_1$, and from
(\ref{ZyChapter19}) in case $\Gamma_S\subseteq S_3.$ Consequently,
for each $\Gamma_S\in {\cal B}(S\setminus S_2)$ and $\Gamma_A\in
{\cal B}(A),$
\begin{eqnarray}
&&M_{\gamma}^{0,\pi}(\Gamma_S\times \Gamma_A)=
\int_{\Gamma_S}M_{\gamma}^{0,\pi}(dx\times
A)\sigma_1(\Gamma_A|x)=\int_{\Gamma_S}\textbf{P}_\gamma^\sigma(X_0\in
dx)
\sigma_1(\Gamma_A|x)\nonumber\\
&=&\textbf{P}_{\gamma}^\sigma(X_0\in \Gamma_S,~A_1\in\Gamma_A),
\end{eqnarray}
where the first equality is by (\ref{ZyChapter12}).

Suppose that (\ref{ZyChapter13}) holds for all $n\le k$. Consider
the case of $n=k+1$ as follows.

Note that for each $n=0,1,\dots,$
\begin{eqnarray*}
M_{\gamma}^{n,\pi}(\hat{S}_1\times
A)=0=M_{\gamma}^{n,\pi}(S_2\times A),
\end{eqnarray*}
where the first equality is by Lemma \ref{ZyChapterLem6}, and the
second equality is by (\ref{ZyChapter14}) and (\ref{ZyChapter15}).
This and the inductive supposition imply
\begin{eqnarray*}
\int_{S_2\times A}\frac{\tilde{q}(\Gamma_S|y,a)}{q_y(a)}
\textbf{P}_{\gamma}^\sigma(X_{k}\in dy,~A_{k+1}\in da)=0=
\int_{S_2 \times
A}\frac{\tilde{q}(\Gamma_S|y,a)}{q_y(a)}M_{\gamma}^{k,\pi}(dy\times
da).
\end{eqnarray*}
Consequently, for each $\Gamma_S\in {\cal B}(S\setminus S_2)$,
\begin{eqnarray}\label{ZyChapter16}
&&\textbf{P}_{\gamma}^\sigma(X_{k+1}\in \Gamma_S) =\int_{S\times
A}\frac{\tilde{q}(\Gamma_S|y,a)}{q_y(a)}\textbf{P}_{\gamma}^\sigma(X_k\in
dy,~A_{k+1}\in da)\nonumber\\
&=& \int_{(S\setminus S_2) \times
A}\frac{\tilde{q}(\Gamma_S|y,a)}{q_y(a)}\textbf{P}_{\gamma}^\sigma(X_k\in
dy,~A_{k+1}\in da)
+\int_{S_2\times A}\frac{\tilde{q}(\Gamma_S|y,a)}{q_y(a)}\textbf{P}_{\gamma}^\sigma(X_{k}\in dy,~A_{k+1}\in da)   \nonumber\\
&=& \int_{(S\setminus S_2) \times
A}\frac{\tilde{q}(\Gamma_S|y,a)}{q_y(a)}M_{\gamma}^{k,\pi}(dy\times
da) +\int_{S_2\times
A}\frac{\tilde{q}(\Gamma_S|y,a)}{q_y(a)}{M}_{\gamma}^{k,\pi}(dy\times
da) \nonumber\\
&=& \int_{S\times
A}\frac{\tilde{q}(\Gamma_S|y,a)}{q_y(a)}{M}_{\gamma}^{k,\pi}(dy\times
da).
\end{eqnarray}
On the other hand, for each $\Gamma_S\in {\cal B}(S\setminus
S_2),$
\begin{eqnarray}\label{ZyChapter18}
&&\int_{S\times
A}\frac{\tilde{q}(\Gamma_S|y,a)}{q_y(a)}{M}_{\gamma}^{k,\pi}(dy\times
da)=\int_{S\times A} \frac{\tilde{q}(\Gamma_S|y,a)}{q_y(a)} q_y(a)
m_{\gamma,k}^\pi(dy\times da)\nonumber\\
&=&\int_{S\times A}
\tilde{q}(\Gamma_S|y,a)m_{\gamma,k}^\pi(dy\times
da)=E_\gamma^\pi\left[\int_{t_{k}}^{t_{k+1}}\int_A
\tilde{q}(\Gamma_S|\xi_t,a)\pi(da|\omega,t)dt\right]\nonumber\\
&=& E_\gamma^\pi\left[\int_0^{\theta_{k+1}}\int_A
\tilde{q}(\Gamma_S|x_k,a)\pi_{k+1}(da|x_0,\theta_1,\dots,x_k,t)dt
\right]\nonumber\\
&=& E_\gamma^\pi\left[ \int_0^{\infty}\int_A
\tilde{q}(\Gamma_S|x_k,a)\pi_{k+1}(da|x_0,\theta_1,\dots,x_k,t)
e^{-\int_0^t \int_A
q_{x_k}(a)\pi_{k+1}(da|x_0,\theta_1,\dots,x_k,s)ds}dt \right]\nonumber\\
&=&E_\gamma^\pi\left[\frac{\int_A
\tilde{q}(\Gamma_S|x_k,a)\pi(da|\omega,t_{k+1})}{\int_A
q_{x_{k}}(a)\pi(da|\omega,t_{k+1})}\right],
\end{eqnarray}
where the first and the third equalities are by
(\ref{ZyChapter17}), whereas the second equality follows from the
fact that if $q_y(a)=0$, then
\begin{eqnarray*}
\tilde{q}(\Gamma_S|y,a)=0=\frac{\tilde{q}(\Gamma_S|y,a)}{q_y(a)}
q_y(a)
\end{eqnarray*}
keeping in mind $\frac{0}{0}=0;$ the similar reasoning justifies
the last equality, too. This together with (\ref{ZyChapter16})
shows
\begin{eqnarray}\label{ZyChapter20}
\textbf{P}_{\gamma}^\sigma(X_{k+1}\in
\Gamma_S)=E_\gamma^\pi\left[\frac{\int_A
\tilde{q}(\Gamma_S|x_k,a)\pi(da|\omega,t_{k+1})}{\int_A
q_{x_{k}}(a)\pi(da|\omega,t_{k+1})}\right]
\end{eqnarray}
for each $\Gamma_S\in {\cal B}(S\setminus S_2).$

Now it holds that
\begin{eqnarray}\label{ZyChapter21}
\textbf{P}_{\gamma}^\sigma(X_{k+1}\in\hat S_1)&=&\int_{(S\setminus
\hat{S}_1) \times A} \tilde{q}(\hat{S}_1|x,a) {m}_{k}^\pi(dx\times
da)
+\int_{\hat{S}_1 \times A} \tilde{q}(\hat{S}_1|x,a)m_{k}^\pi(dx\times da) \nonumber \\
&=&0= M_{\gamma}^{k+1,\pi}(\hat{S}_1\times A),
\end{eqnarray}
where the first equality is by the last to the second equality of
(\ref{ZyChapter18}), whereas the second and the last equalities
are by Lemma \ref{ZyChapterLem6}.

One can see that for each $\Gamma_S\in {\cal B}(S\setminus
(S_2\bigcup \hat{S}_1)),$
\begin{eqnarray*}
&&M_{\gamma}^{k+1,\pi}(\Gamma_S\times
A)=E_\gamma^\pi\left[I\{x_{k+1}\in \Gamma_S\}\right]=E_\gamma^\pi\left[E_\gamma^\pi\left[I\{x_{k+1}\in \Gamma_S\}|x_0,\theta_1,\dots,x_k,\theta_{k+1}\right]\right]\\
&=&E_\gamma^\pi\left[\frac{\int_A
\tilde{q}(\Gamma_S|x_k,a)\pi(da|\omega,t_{k+1})}{\int_A
q_{x_{k}}(a)\pi(da|\omega,t_{k+1})}\right]=\textbf{P}_{\gamma}^\sigma(X_{k+1}\in
\Gamma_S),
\end{eqnarray*}
where the first equality is by the last equality of
(\ref{ZyChapter11}) keeping in mind Lemma \ref{ZyChapterLem5} and
(\ref{ZyChapter19}), and the last equality is by
(\ref{ZyChapter20}). This and (\ref{ZyChapter21}) justify that
\begin{eqnarray*}
M_{\gamma}^{k+1,\pi}(\Gamma_S\times
A)=\textbf{P}_{\gamma}^\sigma(X_{k+1}\in \Gamma_S)
\end{eqnarray*}
for each $\Gamma_S\in {\cal B}(S\setminus S_2).$ Now we see
\begin{eqnarray*}
&&M_{\gamma}^{k+1,\pi}(\Gamma_S\times \Gamma_A)=
\int_{\Gamma_S}M_{\gamma}^{k+1,\pi}(dx\times
A)\sigma_{k+2}(\Gamma_A|x)=\int_{\Gamma_S}\textbf{P}_{\gamma}^\sigma(X_{k+1}\in
dx)
\sigma_{k+2}(\Gamma_A|x)\\
&=&\textbf{P}_{\gamma}^\sigma(X_{k+1}\in \Gamma_S,~A_{k+2}\in
\Gamma_A)
\end{eqnarray*}
for each $\Gamma_S\in {\cal B}(S\setminus S_2)$ and $\Gamma_A\in
{\cal B}(A).$ The statement of the theorem is thus proved by
induction. $\hfill\Box$

\begin{remark}\label{ZyChapterRemark1}
The relation (\ref{ZyChapter13}) in Theorem
\ref{ZyChapterTheorem2} does not hold in general either for
$\Gamma_S\subset S_2,$ or for any given policy $\pi\in\Pi;$ even
under Condition \ref{BookWC1}.
\end{remark}

\section{Occupation measure}\label{GuoZhangSec5}

The objective of this section is to show that restricted on a
measurable subset $\zeta\subseteq S$, the measure
$\eta_\gamma^\pi(dx\times A)$ is $\sigma$-finite; see Theorem
\ref{ZyChapterTheorem1} below, where $\zeta$ is defined by
(\ref{ZyChapter47}), whereas the set $\zeta^c$ is easy to deal
with. After some preliminaries, we do this by adapting the
reasoning of \cite{Dufour:2012}, which is for the occupation
measures for the DTMDP model.

\begin{definition}\label{ForwickChapterDef1}
For each policy $\pi\in \Pi$ for the CTMDP model $(S,A,q,\gamma),$
its (undiscounted) occupation measure $\eta^\pi_\gamma$ is the
measure on ${\cal B}(S\times A)$ given by for each $\Gamma_S\in
{\cal B}(S)$ and $\Gamma_A\in {\cal B}(A)$,
\begin{eqnarray*}
\eta_\gamma^\pi(\Gamma_S\times
\Gamma_A):=E_\gamma^\pi\left[\int_0^\infty I\{\xi_t\in
\Gamma_S\}\pi(\Gamma_A|\omega,t)dt \right].
\end{eqnarray*}
\end{definition}

\begin{lemma}
For each policy $\pi$ for the CTMDP model, its (undiscounted)
occupation measure $\eta^\pi_\gamma$ satisfies the following
relation:
\begin{eqnarray}\label{ZyChapter2}
\int_{\Gamma\times A} q_x(a)\eta^\pi (dx\times
da)+Z^\pi(\Gamma)=\gamma(\Gamma)+\int_{S\times
A}\tilde{q}(\Gamma|y,a)\eta_\gamma^{\pi}(dy\times da)
\end{eqnarray}
for each $\Gamma\in {\cal B}(S),$ where $Z^\pi(\Gamma)\in[0,1].$
\end{lemma}

\par\noindent\textit{Proof.} For each $\alpha>0,$ consider the
measure on ${\cal B}(S\times A)$
\begin{eqnarray}\label{FeinbergChapter13}
\eta_\gamma^{\alpha,\pi}(dx\times da)=
E_\gamma^\pi\left[\int_0^\infty e^{-\alpha t}I\{\xi_t\in
dx\}\pi(da|\omega,t)dt\right],
\end{eqnarray}
which is the ($\alpha$-discounted) occupation measure of the
policy $\pi$ for the CTMDP model $\{S,A,q,\gamma\}$. It follows
from the definition that for each $\pi\in \Pi,$
\begin{eqnarray*}
(\alpha+q_x(a))\eta^{\alpha,\pi}_\gamma(dx\times
da)=\sum_{n=0}^\infty M_{\gamma,\alpha}^{n,\pi}(dx\times da)
\end{eqnarray*}

By Proposition \ref{FeinbergChapterTheorem1}, there is some policy
$\sigma$ for the DTMDP model $\{S_\infty,A,p_\alpha,\gamma\}$
satisfying (\ref{ZyChapter73}) on ${\cal B}(S\times A).$ Note that
$\sum_{n=0}^\infty P_\gamma^{\alpha,\sigma}(X_n\in dx,~A_{n+1}\in
da)$, the right hand side of (\ref{ZyChapter73}), is the
undiscounted occupation measure for the DTMDP model
$\{S_\infty,A,p_\alpha,\gamma\}$ restricted to $S\times A$, so
that, by a well known and easy-to-see fact from the theory of
DTMDPs, for each $\Gamma \in {\cal B}(S)$,
\begin{eqnarray*}
\sum_{n=0}^\infty P_\gamma^{\alpha,\sigma}(X_n\in \Gamma
)=\gamma(\Gamma )+\int_{S\times A}\frac{\tilde{q}(\Gamma
|x,a)}{\alpha+q_x(a)}\sum_{n=0}^\infty
P_\gamma^{\alpha,\sigma}(X_n\in dx,~A_{n+1}\in da).
\end{eqnarray*}
By (\ref{ZyChapter73}), the above can be written as
\begin{eqnarray}\label{FeinbergChapter14}
\int_{\Gamma \times
A}(q_x(a)+\alpha)\eta_\gamma^{\alpha,\pi}(dx\times da)
=\gamma(\Gamma )+ \int_{S\times A} \tilde{q}(\Gamma |x,a)
\eta_\gamma^{\alpha,\pi}(dx\times da).
\end{eqnarray}
Keeping in mind
\begin{eqnarray*}
\int_{\Gamma\times A}\alpha\eta^{\alpha,\pi}_\gamma(dx\times
da)=E^\pi_\gamma\left[\int_0^\infty \alpha e^{-\alpha t}
I\{\xi_t\in \Gamma\}dt\right]\in[0,1]
\end{eqnarray*}
for each $\alpha\in(0,\infty),$ one can legitimately take the
upper limit as $0<\alpha\downarrow 0$ on the both sides of the
above equality to see that $\eta^\pi_\gamma(dx\times da)$
satisfies that for each $\Gamma\in {\cal B}(S)$
\begin{eqnarray*}
&&\int_{\Gamma\times A}q_x(a)\eta^{\pi}_\gamma(dx\times
da)+\limsup_{0<\alpha\downarrow 0} \int_{\Gamma\times
A}\alpha\eta^{\alpha,\pi}_\gamma(dx\times da)\nonumber\\
&=&\gamma(\Gamma)+\int_{S\times A}
\tilde{q}(\Gamma|y,a)\eta^{\pi}_\gamma(dy\times da),
\end{eqnarray*}
where we have used the fact that
$\eta^{\alpha,\pi}_\gamma(dx\times da)\uparrow
\eta^\pi_\gamma(dx\times da)$ setwise as $\alpha\downarrow 0,$ and
the monotone convergence theorem; see Theorem 2.1 of
Hern{\'a}ndez-Lerma and Lasserre \cite{Hernandez-Lerma:2000}. By
putting $Z^\pi(\Gamma)=\limsup_{0<\alpha\downarrow 0}
\int_{\Gamma\times A}\alpha\eta^{\alpha,\pi}_\gamma(dx\times
da)\in[0,1]$, we see that the statement of the lemma holds.
$\hfill\Box$

\begin{remark}
The relation (\ref{FeinbergChapter14}) was established under the
extra conditions imposed on the growth of the transition rates
$q(dy|x,a)$ in \cite{GuoSong:2011,ABPZY:20102}. The relation
(\ref{ZyChapter2}) was established for certain subsets $\Gamma\in
{\cal B}(S)$ in \cite{GuoMZhang:2013}, where the authors imposed
extra conditions and considered the absorbing models, so that the
term $Z^\pi(\Gamma)$ vanishes for all the ``transient'' subsets
$\Gamma.$

\end{remark}

\begin{lemma}\label{ZyChapterLem1}
Let some feasible policy $\pi$ for problem (\ref{ZyChapter1}) with
a finite value be fixed. Suppose that Condition \ref{BookWC1} is
satisfied, and that $\{B_j,~j=1,2,\dots\}\subseteq {\cal B}(S)$ is
an increasing sequence of open sets satisfying
\begin{eqnarray}\label{ZyChapter9}
\int_{B_j\times A} q_x(a)\eta_\gamma^\pi(dx\times da)<\infty.
\end{eqnarray}
Then the following assertions hold.
\par\noindent(a)
There exists a sequence of open sets
$\{E_j,~j=1,2,\dots\}\subseteq {\cal B}(S)$ such that
\begin{eqnarray*}
Y:=\left\{ x\in S: \forall~a\in A,~\tilde{q}(\bigcup_j
B_j|x,a)+\sum_{i=0}^N c_i(x,a)>0\right\}=\bigcup_j E_j,
\end{eqnarray*}
and for all $j=1,2,\dots,$ $ \int_{E_j\times A}
q_x(a)\eta_\gamma^\pi(dx\times da)<\infty. $

\par\noindent(b) There exists a sequence of open sets
$\{\tilde{E}_j,~j=1,2,\dots\}\subseteq {\cal B}(S)$ such that $
Y=\bigcup_j \tilde{E}_j, $ and for all $j=1,2,\dots,$ $
\eta_\gamma^\pi(\tilde{E}_j\times A)<\infty. $
\end{lemma}

\par\noindent\textit{Proof.} (a) Define for each $l=1,2,\dots$ and
$j=1,2,\dots,$
\begin{eqnarray*}
B_j^{(l)}:=\left\{(x,a)\in S\times A:~
\frac{\tilde{q}(B_j|x,a)+\sum_{i=0}^Nc_i(x,a)}{q_x(a)}>\frac{1}{l}\right\}
\end{eqnarray*}
and
\begin{eqnarray*}
C_j^{(l)}:=\left\{ x\in S:~\forall~a\in A,
~\frac{\tilde{q}(B_j|x,a)+\sum_{i=0}^Nc_i(x,a)}{q_x(a)}>\frac{1}{l}\right\}.
\end{eqnarray*}
From (\ref{ZyChapter2}), we see that (\ref{ZyChapter9}) implies
\begin{eqnarray}\label{ZyChapter3}
\int_{S\times A} \tilde{q}(B_j|x,a)\eta_\gamma^\pi(dx\times
da)<\infty.
\end{eqnarray}
Let $ N(q):=\{(x,a)\in S\times A: q_x(a)=0\}. $ Then for each
$j,l=1,2,\dots,$
\begin{eqnarray*}
&&\int_{B_j^{(l)}}q_x(a)\eta_\gamma^\pi(dx,da)=\int_{B_j^{(l)}\bigcap
(N(q)^c)} q_x(a)\eta_\gamma^\pi(dx\times da)\\
&\le& \int_{B_j^{(l)}\bigcap (N(q)^c)}
q_x(a)\frac{\tilde{q}(B_j|x,a)+\sum_{i=0}^Nc_i(x,a)}{q_x(a)} l
\eta_\gamma^{\pi}(dx\times da)\\
&\le & l \left( \int_{S\times
A}\tilde{q}(B_j|x,a)\eta_\gamma^\pi(dx\times da)+\int_{S\times
A}\sum_{i=0}^N c_i(x,a)\eta_\gamma^\pi(dx\times da)
\right)<\infty,
\end{eqnarray*}
where the last inequality follows from (\ref{ZyChapter3}) and the
assumption of the policy $\pi$ being feasible with a finite value.
Since $C_j^{(l)}\times A\subseteq B_j^{(l)}$, it follows that for
each $l,j=1,2,\dots,$
\begin{eqnarray}\label{ZyChapter4}
\int_{C_j^{(l)}\times A} q_x(a)\eta_\gamma^\pi(dx\times da) \le
\int_{B_j^{(l)}} q_x(a)\eta_\gamma^\pi(dx\times da)<\infty.
\end{eqnarray}
Since $B_j$ is open in $S$ for each $j=1,2,\dots,$
$\tilde{q}(B_j|x,a)+\sum_{i=0}^N c_i(x,a)\ge 0$ is lower
semicontinuous in $(x,a)\in S\times A$ according to Lemma
\ref{ForwickChapterLem2}(b). By Lemma \ref{ForwickChapterLem2}(a),
$\frac{\tilde{q}(B_j|x,a)+\sum_{i=0}^N c_i(x,a)}{q_x(a)}$ is lower
semicontinuous in $(x,a)\in S\times A.$ Now referring to Lemma
\ref{ZyChapterLem3}(a), we see that $C_j^{(l)}$ is open for each
$j,l=1,2,\dots.$

Next, let us show
\begin{eqnarray}\label{ZyChapter5}
\bigcup_{j}\bigcup_l C_j^{(l)}=Y=\left\{ x\in S: \forall~a\in
A,~\tilde{q}(\bigcup_j B_j|x,a)+\sum_{i=0}^N c_i(x,a)>0\right\}
\end{eqnarray}
as follows. By Lemma \ref{ZyChapterLem3}(b), for each
$j=1,2,\dots,$
\begin{eqnarray*}
\bigcup_l C_j^{(l)}=\{x\in S:~\forall~a\in
A,~\frac{\tilde{q}(B_j|x,a)+\sum_{i=0}^N c_i(x,a)}{q_x(a)}>0\},
\end{eqnarray*}
so that
\begin{eqnarray*}
\bigcup_j\bigcup_l C_j^{(l)}\subseteq Y.
\end{eqnarray*}
For the opposite direction of the above relation, we argue as
follows. Let some $y\in Y$ be arbitrarily fixed. Since
$\frac{\tilde{q}(B_j|x,a)+\sum_{i=0}^N c_i(x,a)}{q_x(a)}$ is lower
semicontinuous in $(x,a)\in S\times A$ as explained earlier, and
is increasing in $j=1,2,\dots$ keeping in mind that $\{B_j\}$ is
an increasing sequence, one can refer to Lemma \ref{ZyChapterLem4}
for the following interchange of the order of infimum and limit:
\begin{eqnarray*}
&&\lim_{j\rightarrow \infty} \inf_{a\in
A}\left\{\frac{\tilde{q}(B_j|y,a)+\sum_{i=0}^N
c_i(y,a)}{q_y(a)}\right\}=\inf_{a\in A}\lim_{j\rightarrow
\infty}\left\{\frac{\tilde{q}(B_j|y,a)+\sum_{i=0}^N
c_i(y,a)}{q_y(a)}\right\}\\
&=&\inf_{a\in
A}\left\{\frac{\tilde{q}(\bigcup_jB_j|y,a)+\sum_{i=0}^N
c_i(y,a)}{q_y(a)}\right\}>0,
\end{eqnarray*}
where the last inequality follows from the fact that $y\in Y.$
This implies the existence of some $j=1,2,\dots$ such that for
each $a\in A,$ it holds that $\tilde{q}(B_j|y,a)+\sum_{i=0}^N
c_i(y,a)>0,$ i.e.,
 $y\in \bigcup_l C_j^{(l)}\subseteq
 \bigcup_{j}\bigcup_{l}C_j^{(l)}$.
 Since $y\in Y$ is arbitrarily fixed, this verifies
 \begin{eqnarray*}
 Y\subseteq
 \bigcup_{j}\bigcup_{l}C_j^{(l)}.
 \end{eqnarray*}
Hence, (\ref{ZyChapter5}) holds, which in combination with
(\ref{ZyChapter4}), proves the
 statement; remember that $C_j^{(l)}$ is open for each $j,l=1,2,\dots.$

(b) The proof of this part is similar to the one of part (a).
 Instead of $B_j^{(l)}$ and $C_j^{(l)}$, one should now introduce
 for each $j,l=1,2,\dots,$
 \begin{eqnarray*}
 \tilde{B}_j^{(l)}:=\left\{ (x,a)\in S\times A: \tilde{q}(B_j|x,a)+\sum_{i=0}^N c_i(x,a)>\frac{1}{l}\right\}
 \end{eqnarray*}
 and
 \begin{eqnarray*}
 \tilde{C}_j^{(l)}:=\left\{x\in S:~\forall~a\in A, \tilde{q}(B_j|x,a)+\sum_{i=0}^N
 c_i(x,a)>\frac{1}{l}\right\},
\end{eqnarray*}
so that
\begin{eqnarray*}
\eta_\gamma^\pi(\tilde{B}_j^{(l)})\le \int_{\tilde{B}_j^{(l)}} l
\left(\tilde{q}(B_j|y,a)+\sum_{i=0}^N
c_i(y,a)\right)\eta_\gamma^\pi(dy\times da)<\infty.
\end{eqnarray*}
Consequently, $\eta_\gamma^\pi(\tilde{C}_j^{(l)}\times A)<\infty$
for each $j,l=1,2,\dots.$ It is clear now how to proceed the rest
of the reasoning as in the proof of part (a). $\hfill\Box$


\begin{lemma}\label{ZyChapterLem2}
Suppose Condition \ref{BookWC1} is satisfied. Let some feasible
policy $\pi$ for problem (\ref{ZyChapter1}) with a finite value be
fixed. Consider
\begin{eqnarray*}
W:=\bigcup_{j=1}^\infty W_j,
\end{eqnarray*}
where $W_j$ is defined recursively as follows:
\begin{eqnarray*}
W_1:=\left\{x\in S: \forall~a\in A,~\frac{\sum_{i=0}^N
c_i(x,a)}{q_x(a)}>0\right\};
\end{eqnarray*}
and for each $j=1,2,\dots,$
\begin{eqnarray*}
W_{j+1}:=\left\{x\in S: \forall~a\in A,
~\frac{\tilde{q}(\bigcup_{i=1}^j W_i|x,a)+\sum_{i=0}^N
c_i(x,a)}{q_x(a)}>0\right\}.
\end{eqnarray*}
Then for each $j=1,2,\dots,$ $W_j$ is open in $S$, and so is $W$.
Furthermore, $\eta_\gamma^\pi(dx\times A)$, being restricted to
$W\in {\cal B}(S)$, is a $\sigma$-finite measure on ${\cal B}(W);$
in other words, $\eta_\gamma^\pi(dx\times A)$ is $\sigma$-finite
on $W$.
\end{lemma}

\par\noindent\textit{Proof.} First of all, let us show by induction
that for each $m=1,2,\dots,$ $W_m$ is open, and there exists a
sequence of open sets $\{E_j^{(m)}\}\subseteq {\cal B}(S)$ such
that
\begin{eqnarray}\label{ZyChapter67}
W_{m}=\bigcup_j E_j^{(m)}
\end{eqnarray}
and
\begin{eqnarray}\label{ZyChapter68}
\int_{E_j^{(m)}\times A} q_x(a)\eta_\gamma^\pi(dx\times
da)<\infty.
\end{eqnarray}

By Lemma \ref{ZyChapterLem3}(b), $W_1=\bigcup_j E^{(1)}_j,$ where
for each $j=1,2,\dots,$
\begin{eqnarray*}
E^{(1)}_j=\left\{x\in S:~\forall~a\in A, \frac{\sum_{i=0}^N
c_i(x,a)}{q_x(a)}>\frac{1}{j} \right\}
\end{eqnarray*}
is open because $\frac{\sum_{i=0}^N c_i(x,a)}{q_x(a)}$ is lower
semicontinuous in $(x,a)\in S\times A,$ and Lemma
\ref{ZyChapterLem3}(a). Therefore, $W_1$ is open. Moreover,
\begin{eqnarray*}
\int_{E^{(1)}_j\times A} q_x(a)\eta_\gamma^\pi(dx\times da)<\infty
\end{eqnarray*}
because the policy $\pi$ is feasible with a finite value. Note
that $\{E^{(1)}_j,~j=1,2,\dots\}$ is an increasing sequence of
open sets.

Suppose that for each $k\le n,$ $W_k$ is open, and there exists a
sequence of open sets $\{E_j^{(k)},~j=1,2,\dots\}\subseteq {\cal
B}(S)$ such that $ W_{k}=\bigcup_j E_j^{(k)} $ and $
\int_{E_j^{(k)}\times A} q_x(a)\eta_\gamma^\pi(dx\times
da)<\infty. $ Then
\begin{eqnarray*}
W_{n+1}&=&\left\{x\in S:~\forall~a\in
A,~\frac{\tilde{q}(\bigcup_{i=0}^n W_i|x,a)+\sum_{i=0}^N c_i(x,a))}{q_x(a)}>0\right\}\\
&=&\left\{x\in S:~\forall~a\in A,~\frac{\tilde{q}(\bigcup_{i=0}^n
\bigcup_{j} E_j^{(i)}|x,a)+\sum_{i=0}^N
c_i(x,a)}{q_x(a)}>0\right\}.
\end{eqnarray*}
Note that each of the sets $E_j^{(i)},~i=1,2,\dots,n,~j=1,2,\dots$
is open, and $\int_{E_j^{(i)}\times A}
q_x(a)\eta_\gamma^{\pi}(dx\times da)<\infty$ by the inductive
supposition, so that $\bigcup_{i=0}^n \bigcup_{j} E_j^{(i)}$ can
be rewritten as the union of an increasing sequence of open sets
in $S$, each of which is of finite measure with respect to $\int_A
q_x(a)\eta_\gamma^\pi(dx\times da).$ Therefore, one can refer to
Lemma \ref{ZyChapterLem1}(a) for the existence of a sequence of
open sets $\{E^{(n+1)}_j,j=1,2,\dots\}\subseteq{\cal B}(S)$
satisfying
 $W_{n+1}=\bigcup_j E_j^{(n+1)}$ and $\int_{E_j^{(n+1)}\times A} q_x(a)\eta_\gamma^{\pi}(dx\times
da)<\infty.$  Thus, $W_{n+1}$ is open in $S$, and the inductive
argument is completed.

We now prove the statement of the lemma. Let us rewrite
\begin{eqnarray*}
W_1=\left\{x\in S: \forall~a\in A,~\sum_{i=0}^N
c_i(x,a)>0\right\}.
\end{eqnarray*}
By Lemma \ref{ZyChapterLem3}(b),
\begin{eqnarray*}
W_1=\bigcup_{j=1}^\infty \tilde{E}_j^{(1)},
\end{eqnarray*}
where for each $j=1,2,\dots,$
\begin{eqnarray*}
\tilde{E}_j^{(1)}=\left\{x\in S: \forall~a\in A,~\sum_{i=0}^N
c_i(x,a)>\frac{1}{j}\right\},
\end{eqnarray*}
which is open by Lemma \ref{ZyChapterLem3}(a), and satisfies
\begin{eqnarray*}
\eta_\gamma^\pi(\tilde{E}_j^{(1)}\times A)<\infty
\end{eqnarray*}
by the fact that the policy $\pi$ is feasible with a finite value.
For each $m=2,3,\dots,$ by what was established in the beginning
of this proof, (\ref{ZyChapter67}), (\ref{ZyChapter68}) and Lemma
\ref{ZyChapterLem1}(b), which is applicable since $\bigcup_j
E_j^{(m)}$ can be rewritten as the union of an increasing sequence
of open sets each of finite measure with respect to $\int_A
q_x(a)\eta_\gamma^\pi(dx\times da),$ there exists a sequence of
open sets $\{\tilde{E}_j^{(m)},~j=1,2,\dots\}\subseteq {\cal
B}(S)$ such that $ W_m=\bigcup_j \tilde{E}_j^{(m)}$ and
$\eta_\gamma^{\pi}(\tilde{E}_j^{(m)}\times A)<\infty$ for each
$j=1,2,\dots.$ It follows that the statement to be proved holds.
$\hfill\Box$
\bigskip

\begin{definition}
Let us define the set
\begin{eqnarray}\label{ZyChapter47}
\zeta:=\left\{x\in S:~\inf_{\pi\in\Pi_{DM}} E_x^\pi\left[
\int_0^\infty \int_A \sum_{i=0}^N
c_i(\xi_t,a)\pi(da|\omega,t)dt\right]>0\right\}.
\end{eqnarray}
\end{definition}
Here $\Pi_{DM}$ stands for the class of deterministic Markov
policies for the CTMDP model. Under Condition \ref{BookWC1}, one
can refer to Proposition \ref{ForwickChapterTheormeDP} for that
$\zeta$ is a measurable (in fact, open) subset of $S.$

\begin{theorem}\label{ZyChapterTheorem1}
Suppose Condition \ref{BookWC1} is satisfied. Let some feasible
policy $\pi$ for problem (\ref{ZyChapter1}) with a finite value be
fixed. Then $\eta_\gamma^\pi(dx\times A)$ is $\sigma$-finite on
$\zeta.$
\end{theorem}

\par\noindent\textit{Proof.} By Lemma \ref{ZyChapterLem2}, the statement of this theorem would
be proved if we showed
\begin{eqnarray}\label{ZyChapter6}
\zeta \subseteq W,
\end{eqnarray}
where the set $W$ is defined in the statement of Lemma
\ref{ZyChapterLem2}. To this end, let us argue as follows. The
relation (\ref{ZyChapter6}) automatically holds if $W=S.$ Now
consider $W\ne S.$ Let us arbitrarily fix some $x\in W^c=\bigcap_j
W_j^c.$ By the definition of the sets $W_j$ given in the statement
of Lemma \ref{ZyChapterLem2}, for each $j=1,2,\dots,$ there exists
some $a_j\in A$ such that $\tilde{q}(\bigcup_{i=1}^j
W_i|x,a_j)+\sum_{i=0}^N c_i(x,a_j)=0,$ which implies
\begin{eqnarray*}
\lim_{j\rightarrow \infty}\inf_{a\in A} \left\{
\tilde{q}(\bigcup_{i=1}^j W_i|x,a)+\sum_{i=0}^N
c_i(x,a)\right\}=0.
\end{eqnarray*}
By Lemma \ref{ZyChapterLem4},
\begin{eqnarray}\label{ZyChapter7}
\inf_{a\in A}\left\{\tilde{q}(W|x,a)+\sum_{i=0}^N
c_i(x,a)\right\}=\lim_{j\rightarrow \infty}\inf_{a\in A} \left\{
\tilde{q}(\bigcup_{i=1}^j W_i|x,a)+\sum_{i=0}^N
c_i(x,a)\right\}=0.
\end{eqnarray}
Since the set $W$ is open by Lemma \ref{ZyChapterLem2},
$\tilde{q}(W|x,a)+\sum_{i=0}^N c_i(x,a)$ is lower semicontinuous
on $S\times A.$  By Proposition 7.33 of Bertsekas and Shreve
\cite{Bertsekas:1978}, there exists a measurable mapping from $S$
to $A$, which attains the infimum on the left hand side of
(\ref{ZyChapter7}) for each $x\in S$. This mapping gives a
deterministic stationary policy for the CTMDP model $\{S,A,q\}$,
under which, given the initial state $x\in W^c,$ the controlled
process keeps being absorbed at the set $W^c$ without inducing any
cost. This implies that $x\in \zeta^c.$ Since $x\in W^c$ is
arbitrarily fixed, it holds that $W^c\subseteq \zeta^c,$ i.e.,
\begin{eqnarray*}
\zeta\subseteq W,
\end{eqnarray*}
as required. $\hfill\Box$

\begin{definition}
Suppose Condition \ref{BookWC1} is satisfied. Let us fix a
measurable mapping $\psi^\ast$ from $S$ to $A$ such that whenever
$\zeta^c\ne \emptyset,$
\begin{eqnarray}\label{ZyChapter22}
E_x^{\psi^\ast}\left[\int_0^\infty \left(\sum_{i=0}^N
c_i(\xi_t,\psi^\ast(\xi_t))\right)dt\right]=0
\end{eqnarray}
for each $x\in \zeta^c$; and whenever $S_2\ne \emptyset,$
\begin{eqnarray}\label{ZyChapter23}
\psi^\ast(x)=f^\ast(x)
\end{eqnarray}
for each $x\in S_2$, where $f^\ast$ is defined by
(\ref{ZyChapter14}). Such a mapping, or say it interchangeably a
deterministic stationary policy, $\psi^\ast$ exists by Proposition
\ref{ForwickChapterTheormeDP}.
\end{definition}
Note that it necessarily holds that for each $i=0,1,\dots,N,$
\begin{eqnarray}\label{ZyChapter35}
c_i(x,\psi^\ast(x))=0
\end{eqnarray}
for all $x\in \zeta^c,$ whenever $\zeta^c\ne \emptyset.$

The next statement is a direct consequence of Theorem
\ref{ZyChapterTheorem1} and Corollary 7.27.2 of Bertsekas and
Shreve \cite{Bertsekas:1978}.
\begin{corollary}\label{ZyChapterCorollary1}
Suppose Condition \ref{BookWC1} is satisfied, and consider a
feasible policy $\pi\in \Pi$ for problem (\ref{ZyChapter1}) with a
finite value. Then there exists a stationary policy
$\varphi_\pi\in \Pi$ such that
\begin{eqnarray}\label{ZyChapter29}
\eta_\gamma^\pi(dx\times da)=\varphi_\pi(da|x)\eta_\gamma^\pi(dx
\times A)
\end{eqnarray}
on ${\cal B}( \zeta \times A)$, and
\begin{eqnarray}\label{ZyChapter46}
\varphi_\pi(da|x)=\delta_{\psi^\ast(x)}(da),~\forall~ x\in
\zeta^c.
\end{eqnarray}
\end{corollary}

\par\noindent\textit{Proof.} By Theorem \ref{ZyChapterTheorem1},
there is a sequence of disjoint measurable subsets $\{\zeta_n\}$
of $\zeta$ such that $\zeta=\bigcup_{n}\zeta_n$ and for each $n$,
$\eta_\gamma^\pi(\zeta_n\times A)<\infty.$ Now one can refer to
Corollary 7.27.2 of Bertsekas and Shreve \cite{Bertsekas:1978} for
the existence of the stochastic kernels $\varphi_{n,\pi}$ from
$\zeta_n$ to ${\cal B}(A)$ satisfying $\eta_\gamma^\pi(dx\times
da)=\eta_\gamma^\pi(dx\times A)\varphi_{n,\pi}(da|x)$ on ${\cal
B}(\zeta_n\times A)$ for each $n.$ Now the stochastic kernel
$\varphi_\pi$ from $S$ to ${\cal B}(A)$ defined by
\begin{eqnarray*}
\varphi_\pi(da|x)=\sum_{n}\varphi_{n,\pi}(da|x)I\{x\in
\zeta_n\}+\delta_{\psi^\ast(x)}(da)I\{ x\in \zeta^c\}
\end{eqnarray*}
is the required one for the statement. $\hfill\Box$

\begin{definition}
Let us introduce the occupation measure $\textbf{M}_\gamma^\sigma$
of a policy $\sigma$ for the undiscoutned DTMDP model
$\{S_\infty,A,p,\gamma\}$ as a measure on ${\cal B}(S\times A)$
defined by
\begin{eqnarray}\label{ZyChapter32}
{\bf M}^\sigma_\gamma(\Gamma_S\times \Gamma_A):=
\sum_{n=0}^\infty\textbf{P}_\gamma^\sigma(X_n\in
\Gamma_S,~A_{n+1}\in \Gamma_A)
\end{eqnarray}
for each $\Gamma_S\in {\cal B}(S)$ and $\Gamma_A\in {\cal B}(A).$
Here, as before, the transition probability $p$ is defined by
(\ref{ZyChapter64}) and (\ref{ZyChapter65}).
\end{definition}

The next statement is a consequence of Theorem
\ref{ZyChapterTheorem2} and its proof.

\begin{corollary}\label{ZyChapterCorollary2}
Suppose Condition \ref{BookWC1} is satisfied. Consider a feasible
policy $\pi=(\pi_n)\in \Pi$ with a finite value for the CTMDP
problem (\ref{ZyChapter1}) such that
\begin{eqnarray}\label{ZyChapter69}
\pi_{n+1}(da|x_0,\theta_1,\dots,x_n,s)=\delta_{\psi^\ast(x_n)}(da)
\end{eqnarray}
whenever $x_n\in \zeta^c.$ Then there is a Markov policy $\sigma$
for the DTMDP $\{S_\infty, A,p,\gamma\}$ such that
\begin{eqnarray}\label{ZyChapter37}
\sigma_{n+1}(da|X_0,A_1,\dots,X_n)=\delta_{\psi^\ast(X_n)}(da)
\end{eqnarray}
for each $X_n\in \zeta^c$, and
\begin{eqnarray}\label{ZyChapter33}
q_x(a)\eta^{\pi}_\gamma(dx\times da)=
\textbf{M}_\gamma^{\sigma}(dx\times da)
\end{eqnarray}
on ${\cal B}( \zeta \times A).$ Here the mapping $\psi^\ast$ is
the fixed one satisfying (\ref{ZyChapter22}) and
(\ref{ZyChapter23}).
\end{corollary}
\par\noindent\textit{Proof.} Inspecting the proof of Theorem
\ref{ZyChapterTheorem2}, one can see that any Markov policy
$\sigma=(\sigma_n)$ for the DTMDP model $\{S_\infty,A,p\}$ with
$\sigma_{n+1}(da|x)$ satisfying (\ref{ZyChapter60}) and
(\ref{ZyChapter12}) for each $n=0,1,\dots$ fulfils the conditions
of the statement of Theorem \ref{ZyChapterTheorem2}; and there
exists at least one such policy, which we consider now. On ${\cal
B}(\zeta^c\setminus S_2\times A),$ (\ref{ZyChapter12}) reads that
for each $n=0,1,\dots$
\begin{eqnarray*}
&&q_x(a)m_{\gamma,n}^\pi(dx\times da)=\left(\int_A
m_{\gamma,n}^\pi(dx\times db)q_x(b)\right)  \sigma_{n+1}(da|x)\\
&\Leftrightarrow& q_x(a)m_{\gamma,n}^\pi(dx\times
A)\delta_{\psi^\ast(x)}(da)=m_{\gamma,n}^\pi(dx\times
A)q_x(\psi^\ast(x))\sigma_{n+1}(da|x),
\end{eqnarray*}
where the equivalence is by (\ref{ZyChapter69}). Therefore, one
can always put $\sigma_{n+1}(da|x)=\delta_{\psi^\ast(x)}(da)$ for
each $x\in \zeta^c\setminus S_2$ without violating
(\ref{ZyChapter12}). This together with (\ref{ZyChapter60}) shows
that the policy $\sigma$ satisfies (\ref{ZyChapter37}); recall
(\ref{ZyChapter23}). From the discussion in the beginning of this
proof, this policy $\sigma$ satisfies (\ref{ZyChapter13}), by
summing up both sides of which with respect to $n$, we see that
(\ref{ZyChapter33}) is also fulfilled. The corollary is now
proved. $\hfill\Box$
\bigskip

We end this section with the next lemma.
\begin{lemma}\label{ZyChapterLem7}
Suppose Condition \ref{BookWC1} is satisfied. Let some $\sigma$ be
a policy for the DTMDP model $\{S_\infty, A, p, \gamma\}$ such
that
\begin{eqnarray*}
\int_{S\times A} \textbf{M}_\gamma^{\sigma}(dx\times
da)\frac{c_i(x,a)}{q_x(a)}<\infty
\end{eqnarray*}
for each $i=0,1,\dots,N.$ Suppose that there exists a stationary
policy $\sigma^S$ for the DTMDP model $\{S_\infty, A,p,\gamma\}$
satisfying $\textbf{M}^{\sigma}(dx\times
da)=\textbf{M}_\gamma^{\sigma}(dx\times A)\sigma^S(da|x)$ on
${\cal B}(\zeta\times A)$, and
$\sigma^S(da|x)=\delta_{\psi^\ast(x)}(da)$ for each $x\in
\zeta^c.$ Then
\begin{eqnarray*}
\textbf{M}_\gamma^{\sigma^S}(dx\times da)\le
\textbf{M}_\gamma^{\sigma}(dx\times da)
\end{eqnarray*}
on ${\cal B}(\zeta\times A).$
\end{lemma}

\par\noindent\textit{Proof.} According to Proposition
\ref{ForwickChapterTheormeDP}; see especially
(\ref{ForwickChapter16}), and Proposition 9.10 of
\cite{Bertsekas:1978},
\begin{eqnarray*}
\inf_{\pi\in \Pi_{DM}}E_x^\pi\left[\int_0^\infty \int_A
\sum_{i=0}^N
c_i(\xi_t,a)\pi(da|\omega,t)dt\right]=\inf_{\sigma}\textbf{E}_x^\sigma\left[\sum_{n=0}^\infty
\sum_{i=0}^N\frac{c_i(X_n,A_{n+1})}{q_{X_n}(A_{n+1})}\right]
\end{eqnarray*}
for each $x\in S.$ Thus,
\begin{eqnarray*}
\zeta=\left\{x\in S:
\inf_{\sigma}\textbf{E}_x^\sigma\left[\sum_{n=0}^\infty
\sum_{i=0}^N\frac{c_i(X_n,A_{n+1})}{q_{X_n}(A_{n+1})}\right]>0\right\}.
\end{eqnarray*}
Now one can apply Theorem 3.3 of Dufour et al \cite{Dufour:2012}
for the statement. We remark that in \cite{Dufour:2012}, only
nonnegative finitely valued cost functions were considered for the
concerned DTMDP model. A careful inspection of the reasonings
therein reveal that all the cited statements from
\cite{Dufour:2012} in this paper survive when the cost functions
are nonnegative extended real-valued. $\hfill\Box$

\section{Optimality result}\label{GuoZhangSec6}
The main objective of this section is to show the existence of a
stationary optimal policy for the CTMDP problem
(\ref{ZyChapter1}); see Theorem \ref{ZyChapterTheorem5} below. In
the process, we also justify the reduction of the CTMDP problem
(\ref{ZyChapter1}) to the DTMDP problem (\ref{ZyChapter63}) under
Conditions \ref{BookWC1} and \ref{FinitenessCon}; see Remark
\ref{ZyChapterRemark2} below. To this end, we firstly show the
sufficiency of stationary policies for the CTMDP problem
(\ref{ZyChapter1}); see Theorem \ref{ZyChapterTheorem3}.

\begin{theorem}\label{ZyChapterTheorem3}
Suppose Condition \ref{BookWC1} is satisfied, and consider a
feasible policy $\pi$ with a finite value for the CTMDP problem
(\ref{ZyChapter1}) such that for each $n=0,1,\dots,$
\begin{eqnarray}\label{ZyChapter53}
\pi_{n+1}(da|x_0,\theta_1,\dots,x_n,s)=\delta_{\psi^\ast(x_n)}(da)
\end{eqnarray}
whenever $x_n\in \zeta^c.$ Then the stationary policy
$\varphi_\pi$ for the CTMDP problem (\ref{ZyChapter1}) coming from
Corollary \ref{ZyChapterCorollary1} satisfies
\begin{eqnarray}\label{ZyChapter54}
E_\gamma^{\varphi_\pi}\left[\int_0^\infty \int_A
c_i(\xi_t,a)\varphi_{\pi}(da|\xi_t)dt\right]\le
E_\gamma^\pi\left[\int_0^\infty \int_A
c_i(\xi_t,a)\pi(da|\omega,t)dt\right]
\end{eqnarray}
for each $i=0,1,\dots,N.$
\end{theorem}

\par\noindent\textit{Proof.} The proof goes in several steps.

\underline{Step 1.} We show that the stationary policy
$\varphi_\pi$ satisfies that
\begin{eqnarray*}
\varphi_\pi(B(x)|x)<1
\end{eqnarray*}
for almost all $x\in S_1$ with respect to
$\eta^\pi_\gamma(dx\times A),$ where $B(x)$ is given by
(\ref{ZyChapter24}), and $S_1$ is given by (\ref{ZyChapter19})

It suffices to prove the above claim for the case of
\begin{eqnarray}\label{ZyChapter25}
\eta^\pi_\gamma(S_1\times A)>0
\end{eqnarray}
as follows.

Note that
\begin{eqnarray}\label{ZyChapter26}
&&\infty>\int_{S_1\times A} \eta_\gamma^\pi(dx\times da)
\sum_{i=0}^N c_i(x,a)
=\int_{S_1} \eta_\gamma^\pi(dx\times A) \int_A \varphi_\pi(da|x) \sum_{i=0}^N c_i(x,a)\nonumber\\
&=& \sum_{n=0}^\infty \left.E_\gamma^\pi\left[E_\gamma^\pi\left[
\int_{t_n}^{t_{n+1}} \int_{S_1}I\{x_n\in dx\} \int_A
\varphi_\pi(da|x)\sum_{i=0}^N
c_i(x,a)dt\right|x_0,\theta_1,\dots,x_n \right]
\right]\nonumber\\
&=& \sum_{n=0}^\infty E_\gamma^\pi\left[\int_{S_1}I\{x_n\in dx\}
\int_A \varphi_\pi(da|x)\sum_{i=0}^N c_i(x,a) \left.E_\gamma^\pi
\left[ \theta_{n+1} \right|x_0,\theta_1,\dots,x_n\right]\right].
\end{eqnarray}
Suppose for contraction that
\begin{eqnarray}\label{ZyChapter31}
\varphi_\pi(B(x)|x)=1
\end{eqnarray}
on a measurable subset $\Gamma_1\subseteq S_1$ of positive measure
with respect to $\eta_\gamma^\pi(dx\times A)$. It holds that
\begin{eqnarray}\label{ZyChapter27}
\int_A \varphi_\pi(da|x)\sum_{i=0}^N c_i(x,a)\ge \int_{B(x)}
\varphi_\pi(da|x)\sum_{i=0}^N c_i(x,a)>0
\end{eqnarray}
for each $x\in \Gamma_1\subseteq S_1,$ where the last inequality
is by (\ref{ZyChapter24}) and (\ref{ZyChapter19}).

According to (\ref{ZyChapter25}), there exists some $n=0,1,\dots$
such that
\begin{eqnarray*}
P^\pi_\gamma(x_n\in \Gamma_1)>0;
\end{eqnarray*}
and for this $n,$ it must hold that
\begin{eqnarray}\label{ZyChapter28}
E_\gamma^\pi \left[
\theta_{n+1}|x_0,\theta_1,\dots,x_n\right]<\infty
\end{eqnarray}
for almost all $\omega\in\{\omega\in \Omega:x_n(\omega)\in
\Gamma_1\}$ with respect to $P^{\pi}_\gamma(d\omega)$, for
otherwise this together with (\ref{ZyChapter27}) would contradict
the first inequality of (\ref{ZyChapter26}).

The definition of $B(x)$ given by (\ref{ZyChapter24}) and the
inequality (\ref{ZyChapter28}) imply that
\begin{eqnarray}\label{ZyChapter30}
\eta_\gamma^\pi(\{(x,a):x\in \Gamma_1, a\in A\setminus B(x)\})>0,
\end{eqnarray}
where the set in the bracket is measurable because so is the set
$\{(x,a):x\in \Gamma_1, a\in B(x)\}$ according to e.g., Theorem
3.1 of Feinberg et al \cite{FeinbergKZ:2013}. Since
\begin{eqnarray*}
&&\int_{\Gamma_1} \varphi_\pi(A\setminus
B(x)|x)\eta_\gamma^\pi(dx\times A)\\
&=&\int_{\Gamma_1\bigcap \zeta}\varphi_\pi(A\setminus
B(x)|x)\eta_\gamma^\pi(dx\times A)+\int_{\Gamma_1\bigcap
(\zeta^c)}\varphi_\pi(A\setminus B(x)|x)\eta_\gamma^\pi(dx\times
A)\\
&=&\eta_\gamma^\pi(\{(x,a):x\in \Gamma_1, a\in A\setminus B(x)\}),
\end{eqnarray*}
which holds by (\ref{ZyChapter29}) and (\ref{ZyChapter53}), the
relation (\ref{ZyChapter30}) implies that $\varphi_\pi(A\setminus
B(x)|x)>0$  on some measurable subset of
$\Gamma_2\subseteq\Gamma_1$ of positive measure with respect to
$\eta_\gamma^\pi(dx\times A)$. This is a desired contradiction
against the relation in (\ref{ZyChapter31}). Step 1 is completed.

\underline{Step 2.} Consider the policy $\sigma$ for the DTMDP
model $\{S_\infty, A, p, \gamma\}$ from Corollary
\ref{ZyChapterCorollary2}, and define the stationary policy
$\sigma^S$ for the DTMDP model $\{S_\infty,A,p,\gamma\}$ by
\begin{eqnarray*}
\sigma^S(da|x)=\delta_{\psi^\ast(x)}(da)
\end{eqnarray*}
for all $x\in \zeta^c$ and for all $x\in \zeta$ satisfying $\int_A
q_x(a)\varphi_\pi(da|x)=0$; and
\begin{eqnarray}\label{ZyChapter34}
\sigma^S(da|x):= \frac{q_x(a)\varphi_\pi(da|x)}{\int_A
q_x(a)\varphi_\pi(da|x)}
\end{eqnarray}
for all $x\in \zeta$ such that $\int_A q_x(a)\varphi_\pi(da|x)>0$.
Recall that $\psi^\ast$ is the fixed measurable mapping satisfying
(\ref{ZyChapter22}) and (\ref{ZyChapter23}). We verify that
\begin{eqnarray*}
\textbf{M}_\gamma^{\sigma}(dx\times A)\sigma^S(da|x)=
\textbf{M}_\gamma^{\sigma}(dx\times da)
\end{eqnarray*}
on ${\cal B}(\zeta\times A);$ recall (\ref{ZyChapter32}) for the
definition of $\textbf{M}_\gamma^\sigma.$ Throughout the proof of
this theorem, the policies $\sigma$ and $\sigma^S$ are understood
as here.

Indeed, on ${\cal B}(\zeta\times A),$ it holds that
\begin{eqnarray*}
&&\textbf{M}_\gamma^{\sigma}(dx\times A)\sigma^S(da|x)=
\left(\int_A \eta_\gamma^\pi(dx\times db) q_x(b)\right)\sigma^S(da|x)\\
&=& \left(\int_{A} \eta_\gamma^{\pi}(dx\times A)\varphi_\pi(db|x)
q_x(b)\right)\sigma^S(da|x)\\
&=& \eta_\gamma^{\pi}(dx\times A) \left(\int_A \varphi_\pi(db|x)
q_x(b)\right) \frac{q_x(a)\varphi_\pi(da|x)}{\int_A
q_x(a)\varphi_\pi(da|x)}\\
&=& \eta_\gamma^{\pi}(dx\times A)q_x(a)\varphi_\pi(da|x)\\
&=&\textbf{M}_\gamma^{\sigma}(dx\times da),
\end{eqnarray*}
where the first and the last equalities are by
(\ref{ZyChapter33}), the second equality is by
(\ref{ZyChapter29}), the third and forth equalities are by
(\ref{ZyChapter34}); and the fact that $\int_A
q_x(a)\varphi_\pi(da|x)>0$ for almost all $x\in \zeta$, which in
turn follows from the facts that $\int_A q_x(a)\psi_\pi(da|x)>0$
for almost all $x\in S_1$ with respect to
$\eta_\gamma^\pi(dx\times A)$ as established in Step 1; $\int_A
q_x(a)\psi_\pi(da|x)>0$ for all $x\in S_3$ by (\ref{ZyChapter19});
and the relation $\zeta\subseteq S_1\bigcup S_3.$ Step 2 is thus
completed.

\underline{Step 3.} We verify that
\begin{eqnarray}\label{ZyChapter36}
\textbf{M}_\gamma^{\sigma^S}(dx\times
A)\sigma^S(da|x)=\textbf{M}_\gamma^{\sigma^S}(dx\times da)\le
\textbf{M}_\gamma^\sigma(dx\times da)
\end{eqnarray}
on ${\cal B}(\zeta\times A)$ as follows.

The equality in (\ref{ZyChapter36}) holds because the policy
$\sigma^S$ is stationary and (\ref{ZyChapter32}). For the
inequality in (\ref{ZyChapter36}), we observe that
\begin{eqnarray*}
&&\int_{S\times A} \textbf{M}_\gamma^{\sigma}(dx\times da)
\frac{c_i(x,a)}{q_x(a)}\\
&=& \int_{\zeta\times A} \textbf{M}_\gamma^{\sigma} (dx\times da)
\frac{c_i(x,a)}{q_x(a)}+\int_{\zeta^c} \textbf{M}_\gamma^{\sigma}
(dx\times A)
\frac{c_i(x,\psi^\ast(x))}{q_x(\psi^\ast(x))}\\
&=& \int_{\zeta\times A}
 \eta_\gamma^{\pi} (dx\times da)
q_x(a)\frac{c_i(x,a)}{q_x(a)}\le \int_{\zeta\times A}\eta_\gamma^\pi(dx\times da)c_i(x,a)\\
&\le& \int_{S\times A} c_i(x,a)\eta_\gamma^{\pi}(dx\times
da)<\infty
\end{eqnarray*}
for each $i=0,1,\dots,N$, where the first equality is by
(\ref{ZyChapter37}), the second equality is by
(\ref{ZyChapter35}), the first inequality is by that
$q_x(a)\frac{c_i(x,a)}{q_x(a)}\le c_i(x,a)$; recall the convention
of $\frac{0}{0}=0$ and $0\cdot\infty=0$, and the last inequality
is by that the policy $\pi$ is feasible with a finite value for
problem (\ref{ZyChapter1}). With this inequality and the equality
of (\ref{ZyChapter36}) in hand, we see that the conditions of
Lemma \ref{ZyChapterLem7} are satisfied, following from which, the
inequality of (\ref{ZyChapter36}) holds. Step 3 is completed.

\underline{Step 4.} Let us introduce the set
\begin{eqnarray}\label{ZyChapter70}
\zeta_\pi:=\left\{x\in \zeta: \int_A q_x(a)\varphi_\pi(da|x)=0
\right\},
\end{eqnarray}
which is measurable. We establish
\begin{eqnarray}\label{ZyChapter42}
\eta_\gamma^{\varphi_\pi}(dx\times
da)q_x(a)=\textbf{M}_\gamma^{\sigma^S}(dx\times da)
\end{eqnarray}
on ${\cal B}((\zeta\setminus\zeta_\pi)\times A).$

To this end, we show by induction the more detailed relation
\begin{eqnarray}\label{ZyChapter40}
M_\gamma^{n,\varphi_\pi}(dx\times da)=
\textbf{P}_\gamma^{\sigma^S}(X_n\in dx,A_{n+1}\in da)
\end{eqnarray}
on ${\cal B}(\zeta\setminus \zeta_\pi\times A)$ for each
$n=0,1,\dots$ as follows.

Consider $n=0$. Then on ${\cal B}(\zeta\setminus \zeta_\pi)$,
\begin{eqnarray}\label{ZyChapter39}
M_\gamma^{0,\varphi_\pi}(dx\times
A)=\gamma(dx)=\textbf{P}_\gamma^{\sigma^S}(X_0\in dx),
\end{eqnarray}
where the first equality is by (\ref{ZyChapter38}) and the fact
that $\zeta\subseteq S\setminus S_2.$ Now on ${\cal
B}(\zeta\setminus \zeta_\pi\times A),$
\begin{eqnarray}\label{ZyChapter41}
&&\textbf{P}_\gamma^{\sigma^S}(X_0\in dx, A_1\in da)=
\textbf{P}_\gamma^{\sigma^S}(X_0\in dx)\sigma^S(da|x)= M_\gamma^{0,\varphi_\pi}(dx\times A)\sigma^S(da|x)\nonumber\\
&=&\int_A E^{\varphi_\pi}_\gamma\left[\int_0^{t_1}
q_x(b)\varphi_\pi(db|x)I\{x_0\in dx\} dt
\right]\frac{q_x(a)\varphi_\pi(da|x)}{\int_A
q_x(a)\varphi_\pi(da|x)}\nonumber\\
&=& E^{\varphi_\pi}_\gamma\left[\int_0^{t_1}
q_x(a)\varphi_\pi(da|x)I\{x_0\in dx\} dt \right]\nonumber\\
&=&M_\gamma^{0,\varphi_\pi}(dx\times da),
\end{eqnarray}
where the second equality is by (\ref{ZyChapter39}), the third
equality is by (\ref{ZyChapter34}); remember that
\begin{eqnarray*}
\int_A q_x(a)\varphi_\pi(da|x)>0,~\forall~x\in \zeta\setminus
\zeta_\pi.
\end{eqnarray*}

Assume (\ref{ZyChapter40}) holds on ${\cal B}(\zeta\setminus
\zeta_\pi\times A)$ for all $n\le k$, and consider the case of
$n=k+1.$ On the one hand, on ${\cal B}(\zeta\setminus \zeta_\pi)$
it holds that
\begin{eqnarray*}
&&\textbf{P}_\gamma^{\sigma^S}(X_{k+1}\in dx) =\int_{S\times A}
\frac{\tilde{q}(dx|y,a)}{q_y(a)}
\textbf{P}^{\sigma^S}_\gamma(X_k \in dy,~A_{k+1}\in da)\\
&=&\int_{S\times A} \frac{\tilde{q}(dx|y,a)}{q_y(a)}
M_\gamma^{k,\varphi_\pi}(dy\times da)=\int_{S\times A}
\frac{\tilde{q}(dx|y,a)}{q_y(a)}q_y(a)
m_{\gamma,k}^{\varphi_\pi}(dy\times da)\\
&=& \int_{S\times A} \tilde{q}(dx|y,a)
m_{\gamma,k}^{\varphi_\pi}(dy\times
da)\\
&=&E^{\varphi_\pi}_\gamma\left[\int_A
\tilde{q}(dx|x_{k},a)\varphi_\pi(da|x_{k})
\left.E_\gamma^\varphi\left[\theta_{k+1} \right|x_0,\theta_1,\dots,x_k\right]\right]\\
&=&E^\varphi_\gamma\left[\frac{\int_A
\tilde{q}(dx|x_{k},a)\varphi_\pi(da|x_{k})}{\int_A
 q_{x_{k}}(a)\varphi_\pi(da|x_{k}) }\right],
\end{eqnarray*}
where the second equality is by the inductive supposition, the
forth equality is by that
$\frac{\tilde{q}(dx|y,a)}{q_y(a)}q_y(a)=\tilde{q}(dx|y,a)$ no
matter whether $q_y(a)$ vanishes or not, and the last equality
holds due to the convention of $\frac{0}{0}=0.$ On the other hand,
on ${\cal B}(\zeta\setminus \zeta_\pi)$,
\begin{eqnarray*}
&&M_{k+1,\gamma}^{\varphi_\pi}(dx\times A)\\
&=&\left.E_\gamma^{\varphi_\pi}\left[E_\gamma^{\varphi_\pi}\left[I\{x_{k+1}\in
dx\}\right.\right.\right.\\
&&\left.\left.\left.\left.E_\gamma^{\varphi_\pi}\left[\int_0^{\theta_{k+2}}\int_A
q_{x_{k+1}}(a)\varphi_\pi(da|x_{k+1}) \right|x_0,\theta_1,\dots,x_k,\theta_{k+1},x_{k+1}\right] \right|x_0,\theta_1,\dots,x_k,\theta_{k+1}\right]\right]\\
&=&\left.E_\gamma^{\varphi_\pi}\left[E_\gamma^{\varphi_\pi}\left[
I\{x_{k+1}\in
dx\}\right|x_0,\theta_1,\dots,x_k,\theta_{k+1} \right]\right]\\
&=&E^{\varphi_\pi}_\gamma\left[\frac{\int_A
\tilde{q}(dx|x_{k},a)\varphi_\pi(da|x_{k})}{\int_A
 (q_{x_{k}}(a))\varphi_\pi(da|x_{k}) }\right].
\end{eqnarray*}
Thus,
\begin{eqnarray*}
M_{k+1,\gamma}^{\varphi_\pi}(dx\times A)=
\textbf{P}_\gamma^{\sigma^S}(X_{k+1}\in dx)
\end{eqnarray*}
on ${\cal B}(\zeta\setminus \zeta_\pi)$. Based on this, a similar
calculation as the one for (\ref{ZyChapter41}) leads to
\begin{eqnarray*}
M_{k+1,\gamma}^{\varphi_\pi}(dx\times da)=
\textbf{P}_\gamma^{\sigma^S}(X_{k+1}\in dx,~A_{k+1}\in da)
\end{eqnarray*}
on ${\cal B}(\zeta\setminus\zeta_\pi)\times A).$ Hence,
(\ref{ZyChapter40}) is shown by induction, and (\ref{ZyChapter42})
follows. Step 4 is completed.

\underline{Step 5.} We show that
\begin{eqnarray}\label{ZyChapter44}
\eta_\gamma^{\varphi_\pi}(dx\times da)\le\eta_\gamma^\pi(dx\times
da)
\end{eqnarray}
on ${\cal B}(\zeta\setminus\zeta_\pi\times A).$

Indeed, by (\ref{ZyChapter36}) and (\ref{ZyChapter42}) as
established in Steps 3 and 4, we see
\begin{eqnarray*}
\eta_\gamma^{\varphi_\pi}(dx\times da)q_x(a)\le
\textbf{M}_\gamma^\sigma(dx\times da)
\end{eqnarray*}
on ${\cal B}(\zeta\setminus \zeta_\pi\times A)$, which together
with (\ref{ZyChapter33}) further leads to
\begin{eqnarray}\label{ZyChapter43}
\eta_\gamma^{\varphi_\pi}(dx\times da)q_x(a)\le
\eta_\gamma^\pi(dx\times da)q_x(a)
\end{eqnarray}
on ${\cal B}(\zeta\setminus \zeta_\pi\times A)$. Now on ${\cal
B}(\zeta\setminus\zeta_\pi)$,
\begin{eqnarray*}
&&E_\gamma^{\varphi_\pi}\left[\int_0^\infty \int_A
q_x(a)\varphi_\pi(da|x)I\{\xi_t\in dx\}dt \right]\\
&=&\left(\int_A
\varphi_\pi(da|x)q_x(a)\right)\eta_\gamma^{\varphi_\pi}(dx\times
A)=\int_A \eta_\gamma^{\varphi_\pi}(dx\times da)q_x(a) \le \int_A
\eta_\gamma^{\pi}(dx\times da)q_x(a)\\
&=&\left(\int_A
q_x(a)\varphi_\pi(da|x)\right)\eta_\gamma^\pi(dx\times A),
\end{eqnarray*}
where the inequality is by (\ref{ZyChapter43}), and the last
equality is by (\ref{ZyChapter29}). Since $\int_A
q_x(a)\varphi_\pi(da|x)>0$ for all $x\in \zeta\setminus
\zeta_\pi,$ we infer from the above inequality for that
\begin{eqnarray*}
\eta_\gamma^{\varphi_\pi}(dx\times A)\le \eta_\gamma^\pi(dx\times
A)
\end{eqnarray*}
on ${\cal B}(\zeta\setminus\zeta_\pi),$ from which
(\ref{ZyChapter44}) holds on ${\cal
B}(\zeta\setminus\zeta_\pi\times A);$ recall (\ref{ZyChapter29}).
Step 5 is completed.

\underline{Step 6.}  We show that
\begin{eqnarray}\label{ZyChapter52}
\eta_\gamma^{\varphi_\pi}(\zeta_\pi\times A)=0.
\end{eqnarray}

Suppose for contradiction that
\begin{eqnarray}\label{ZyChapter45}
\eta_\gamma^{\varphi_\pi}(\zeta_\pi\times A)>0.
\end{eqnarray}
Note that $\zeta_\pi\subseteq S_1,$ where $\zeta_\pi$ is given by
(\ref{ZyChapter70}); recall that $\zeta\subseteq S_1\bigcup S_3$
and the definition of $S_3.$ Therefore, the statement established
in Step 1 implies that
\begin{eqnarray}\label{ZyChapter51}
\eta_\gamma^\pi(\zeta_\pi\times A)=0.
\end{eqnarray}
Therefore, $\gamma(\zeta_\pi)=0.$ Now following from
(\ref{ZyChapter45}), there exists some $\Gamma\in {\cal
B}(S\setminus \zeta_\pi)$ satisfying that
\begin{eqnarray}\label{ZyChapter48}
\int_A \tilde{q}(\zeta_\pi|x,a)\varphi_\pi(da|x)>0
\end{eqnarray}
for all $x\in \Gamma,$ and
\begin{eqnarray}\label{ZyChapter49}
\eta_\gamma^{\varphi_\pi}(\Gamma\times A)>0.
\end{eqnarray}
Note that according to (\ref{ZyChapter46}), the definition of the
set $\zeta$ given by (\ref{ZyChapter47}), and (\ref{ZyChapter22}),
we see that $\tilde{q}(\zeta|x,\psi^\ast(x))=0$ for each $x\in
\zeta^c$. Since $\zeta_\pi\subseteq \zeta$, we see
$\tilde{q}(\zeta_\pi|x,\psi^\ast(x))=0$ for each $x\in \zeta^c$.
Consequently, we have
\begin{eqnarray*}
\Gamma\in {\cal B}(\zeta\setminus \zeta_\pi)
\end{eqnarray*}
for otherwise it would contradict (\ref{ZyChapter48}). This fact,
(\ref{ZyChapter49}) and (\ref{ZyChapter44}) as established in Step
5 show that
\begin{eqnarray}\label{ZyChapter50}
\eta_\gamma^\pi(\Gamma\times A)>0.
\end{eqnarray}
Now
\begin{eqnarray*}
\int_{\Gamma\times A}\eta_\gamma^\pi(dx\times
da)\tilde{q}(\zeta_\pi|x,a)= \int_\Gamma \eta_\gamma^\pi(dx\times
A)\int_A \tilde{q}(\zeta_\pi|x,a)\varphi_\pi(da|x)>0
\end{eqnarray*}
where the first equality is by (\ref{ZyChapter29}), and the last
inequality is by (\ref{ZyChapter48}) and (\ref{ZyChapter50}).
Thus,
\begin{eqnarray*}
E_\gamma^\pi\left[\int_0^\infty \int_A
\tilde{q}(\zeta_\pi|\xi_t,a)\pi(da|\omega,t)I\{\xi_t\in \Gamma\}dt
\right]>0.
\end{eqnarray*}
It follows from this inequality and the construction of the CTMDP
that $\eta_\gamma^\pi(\zeta_\pi\times A)>0,$ which is a
contradiction against (\ref{ZyChapter51}). Hence,
(\ref{ZyChapter52}) holds. Step 6 is completed.

\underline{Step 7.} We prove the statement of the theorem now. It
holds that for each $i=0,1,\dots,N,$
\begin{eqnarray*}
&&\int_{S\times A} \eta_\gamma^\pi(dx\times da)
c_i(x,a)\\
&=& \int_{\zeta\setminus \zeta_\pi\times A}
\eta_\gamma^\pi(dx\times da) c_i(x,a)+\int_{\zeta^c\times A}
\eta_\gamma^\pi(dx\times da) c_i(x,\psi^\ast(x))
+\int_{\zeta_\pi\times A}
\eta_\gamma^\pi(dx\times da) c_i(x,a)\\
&\ge& \int_{\zeta\setminus \zeta_\pi\times A}
\eta_\gamma^{\varphi_\pi}(dx\times da)
c_i(x,a)+\int_{\zeta^c\times A} \eta_\gamma^{\varphi_\pi}(dx\times
da) c_i(x,\psi^\ast(x)) +\int_{\zeta_\pi\times A}
\eta_\gamma^{\varphi_\pi}(dx\times da) c_i(x,a)\\
&=& \int_{S\times A} \eta_\gamma^{\varphi_\pi}(dx\times da)
c_i(x,a),
\end{eqnarray*}
where the first equality is by (\ref{ZyChapter53}), and the
inequality is by (\ref{ZyChapter35}), (\ref{ZyChapter46}),
(\ref{ZyChapter44}), and (\ref{ZyChapter52}). Thus,
(\ref{ZyChapter54}) is proved. $\hfill\Box$

\begin{corollary}\label{ZyChapterCorollary3}
Suppose Condition \ref{BookWC1} is satisfied, and consider a
feasible policy $\pi$ with a finite value for the CTMDP problem
(\ref{ZyChapter1}) satisfying (\ref{ZyChapter53}) as in the
statement of Theorem \ref{ZyChapterTheorem3}. Then there exists a
stationary policy $\phi_\pi$ such that
\begin{eqnarray}\label{ZyChapter57}
\phi_\pi(B(x)|x)=0
\end{eqnarray}
for each $x\in S_1\setminus \hat{S}_1$ provided that $S_1\setminus
\hat{S}_1\ne \emptyset,$
\begin{eqnarray}\label{ZyChapter59}
\phi_\pi(da|x)=\delta_{\psi^\ast(x)}(da)
\end{eqnarray}
for each $x\in \zeta^c$ whenever $\zeta^c\ne \emptyset$, and
\begin{eqnarray*}
E_\gamma^{\phi_\pi}\left[\int_0^\infty \int_A
c_i(\xi_t,a)\phi_{\pi}(da|\xi_t)dt\right]\le
E_\gamma^\pi\left[\int_0^\infty \int_A
c_i(\xi_t,a)\pi(da|\omega,t)dt\right]
\end{eqnarray*}
for each $i=0,1,\dots,N.$
\end{corollary}

\par\noindent\textit{Proof.}
Let the stationary policy $\varphi_\pi$ be as in the statement of
Theorem \ref{ZyChapterTheorem3}. For each $x\in S_1\setminus
\hat{S}_1$, $A\setminus B(x)\ne \emptyset$; this is by the
definitions of $B(x),$ $S_1$ and $\hat{S}_1$; see
(\ref{ZyChapter24}) and (\ref{ZyChapter19}). By Proposition 7.33
of Bertsekas and Shreve \cite{Bertsekas:1978}, there is a
measurable mapping $\hat{\psi}$ from $S_1\setminus \hat{S}_1$ to
$A$ such that
\begin{eqnarray*}
\sup_{a\in A}q_x(a)=q_x(\hat{\psi}(x))>0
\end{eqnarray*}
for each $x\in S_1\setminus \hat{S}_1,$ where the inequality
follows from the fact that $\sup_{a\in A}q_x(a)=\max_{a\in
A}q_x(a)=\max_{a\in A\setminus B(x)}q_x(a)>0$; recall the
definition of $B(x)$ as given by (\ref{ZyChapter24}). Observe that
\begin{eqnarray}\label{ZyChapter71}
\{x\in S_1\setminus \hat{S}_1:\varphi_\pi(A\setminus
B(x)|x)=0\}=\{x\in (S_1\setminus \hat{S}_1)\bigcap
\zeta:\varphi_\pi(A\setminus B(x)|x)=0\}
\end{eqnarray}
by (\ref{ZyChapter22}) and the definition of $S_1.$ Now if
\begin{eqnarray*}
\{x\in S_1\setminus \hat{S}_1:\varphi_\pi(A\setminus
B(x)|x)=0\}\ne \emptyset,
\end{eqnarray*}
then we modify the definition of $\varphi_\pi$ by putting (with
slight abuse of notations by using $\varphi_\pi$ for both the
original and the modified policies)
$\varphi_\pi(da|x):=\delta_{\hat{\psi}(x)}(da)$ for each $x\in
\{x\in (S_1\setminus \hat{S}_1)\bigcap
\zeta:\varphi_\pi(A\setminus B(x)|x,a)=0\}.$ Since
\begin{eqnarray*}
\eta_\gamma^\pi(\{x\in S_1\setminus
\hat{S}_1:\varphi_\pi(A\setminus B(x)|x)=0\})=0
\end{eqnarray*}
as established in Step 1 of the proof of Theorem
\ref{ZyChapterTheorem3}, the resulting stationary policy
$\varphi_\pi$ still satisfies (\ref{ZyChapter29}) and
(\ref{ZyChapter46}); recall (\ref{ZyChapter71}). Therefore,
Theorem \ref{ZyChapterTheorem3} remains applicable to this
modified policy. For this reason, in the rest of this proof, we
suppose without loss of generality that
\begin{eqnarray}\label{ZyChapter72}
\{x\in S_1\setminus \hat{S}_1:\varphi_\pi(A\setminus B(x)|x)=0\}=
\emptyset.
\end{eqnarray}

Now define a stationary policy $\phi_\pi$ by
\begin{eqnarray*}
\phi_\pi(da|x):=\frac{\varphi_\pi(da\bigcap (A\setminus
B(x))|x)}{\varphi_\pi( (A\setminus B(x))|x)}
\end{eqnarray*}
for each $x\in S_1\setminus \hat{S}_1$, and
\begin{eqnarray*}
\phi_\pi(da|x):=\varphi_\pi(da|x)
\end{eqnarray*}
elsewhere. Observe that $\phi_\pi$ defined in the above is indeed
a stochastic kernel; this follows from the fact that
$\{(x,a):q_x(a)=0\}=\{(x,a):a\in B(x)\}$ is measurable, which is
by Theorem 3.1 of Feinberg et al \cite{FeinbergKZ:2013}; see also
Corollary 18.8 of \cite{Aliprantis:2007}, and Proposition 7.29 of
\cite{Bertsekas:1978}. The relation (\ref{ZyChapter59}) holds for
this policy $\phi_\pi$ because of its definition and
(\ref{ZyChapter46}); observe that for each $x\in (S_1\setminus
\hat{S}_1)\bigcap (\zeta^c),$ it holds that $\psi^\ast(x)\notin
B(x).$

Direct calculations show that for each $x\in S,$
\begin{eqnarray*}
\frac{\int_A\tilde{q}(dy|x,a)\phi_\pi(da|x) }{\int_A q_x(a)
\phi_\pi(da|x)}=\frac{\int_A\tilde{q}(dy|x,a)\varphi_\pi(da|x)
}{\int_A q_x(a) \varphi_\pi(da|x)}.
\end{eqnarray*}
Also observe that for each $x\in S_1\setminus \hat{S}_1$ and
$i=0,1,\dots,N,$
\begin{eqnarray*}
&&\int_0^\infty \int_{A} c_i(x,a)\phi_\pi(da|x) e^{-
\int_{A}q_x(a)\phi_\pi(da|x)t} dt\\
&=&  \int_0^\infty \frac{\int_{A} c_i(x,a)\varphi_\pi(da|x)-
\int_{B(x)} c_i(x,a)\varphi_\pi(da|x) }{\varphi_\pi(A\setminus
B(x)|x)}e^{- \int_{A}q_x(a){\varphi_\pi}(da|x)t
\frac{1}{\varphi_\pi(A\setminus
B(x)|x)}}dt\\
&\le &\int_{A} c_i(x,a){\varphi_\pi}(da|x)
\frac{1}{\int_{A}q_x(a)\varphi_\pi(da|x)}\\
&\le & \int_0^\infty \int_{A} c_i(x,a)\varphi_\pi(da|x) e^{-
\int_{A}q_x(a)\varphi_\pi(da|x)t} dt;
\end{eqnarray*}
remember, $\int_{A}q_x(a){\varphi_\pi}(da|x)>0$ for each $x\in
S_1\setminus \hat{S}_1$ by (\ref{ZyChapter72}). In other words,
under the stationary policy $\phi_\pi,$ given the current state
$x\in S$, the (conditional) distribution of the next jump-in state
is the same as the one under the stationary policy $\varphi_\pi,$
and the total (conditional) expected cost during the current
sojourn time is not larger than the one under $\varphi_\pi.$ Since
both policies $\varphi_\pi$ and $\phi_\pi$ are stationary, this
and Theorem \ref{ZyChapterTheorem3} prove the statement.
$\hfill\Box$

\begin{corollary}\label{ZyChapterCorollary4}
Suppose Condition \ref{BookWC1} is satisfied, and consider a
feasible policy $\pi$ with a finite value for the CTMDP problem
(\ref{ZyChapter1}) satisfying (\ref{ZyChapter53}) as in the
statement of Theorem \ref{ZyChapterTheorem3}. Then there exists a
stationary policy $\sigma^S_\pi$ for the DTMDP model
$\{S_\infty,A,p,\gamma\}$ such that for each $i=0,1,\dots,N,$
\begin{eqnarray*}
\textbf{E}_\gamma^{\sigma^S_\pi}\left[\sum_{n=0}^\infty
\frac{c_i(X_n,A_{n+1})}{q_{X_n}(A_{n+1})}\right]\le
E_\gamma^\pi\left[\int_0^\infty \int_A
c_i(\xi_t,a)\pi(da|\omega,t)dt\right].
\end{eqnarray*}
\end{corollary}

\par\noindent\textit{Proof.} Let $\phi_\pi$ be the stationary
policy for the CTMDP model coming from Corollary
\ref{ZyChapterCorollary3}. By Theorem \ref{ZyChapterTheorem2},
there is a Markov policy say $\sigma^M_\pi=({\sigma^M_\pi}_{n})$
for the DTMDP model $\{S_\infty,A,p,\gamma\}$ satisfying, for each
$n=0,1,\dots,$
\begin{eqnarray}\label{ZyChapter55}
M_\gamma^{n,\phi_\pi}(dx\times
da)=\textbf{P}_\gamma^{\sigma^M_\pi}(X_n\in dx,~A_{n+1}\in da)
\end{eqnarray}
on ${\cal B}(S\setminus S_2\times A)$, and
\begin{eqnarray}\label{ZyChapter61}
{\sigma^M_\pi}_{n+1}(da|x)=\delta_{f^\ast(x)}(da)
\end{eqnarray} for each
$x\in S_2$ whenever $S_2\ne \emptyset.$

Now for each $i=0,1,\dots,N,$ it holds that
\begin{eqnarray}\label{ZyChapter56}
&&E_\gamma^{\phi_\pi}\left[\int_0^\infty \int_A
c_i(\xi_t,a)\phi_\pi(da|\xi_t)dt\right]=\sum_{n=0}^\infty
\int_{S\times A} c_i(x,a)m_{\gamma,n}^{\phi_\pi}(dx\times da)\nonumber \\
&=&\sum_{n=0}^\infty\left\{\int_{S_1\setminus \hat{S}_1\times A}
c_i(x,a)m_{\gamma,n}^{\phi_\pi}(dx\times da)+\int_{\hat{S_1}\times
A} c_i(x,a)m_{\gamma,n}^{\phi_\pi}(dx\times da) \right.\nonumber\\
&&\left.+\int_{S_2\times A}
c_i(x,a)m_{\gamma,n}^{\phi_\pi}(dx\times da)+\int_{S_3\times A}
c_i(x,a)m_{\gamma,n}^{\phi_\pi}(dx\times da)\right\}.
\end{eqnarray}
The first term in the summand in the last line of the above
equality can be written as follows:
\begin{eqnarray*}
&&\int_{S_1\setminus \hat{S}_1\times A}
c_i(x,a)m_{\gamma,n}^{\phi_\pi}(dx\times da)=\int_{S_1\setminus
\hat{S}_1 } \int_A
c_i(x,a)\phi_\pi(da|x)m_{\gamma,n}^{\phi_\pi}(dx\times A)\\
&=&\int_{S_1\setminus \hat{S}_1 } \int_A
\frac{c_i(x,a)}{q_x(a)}q_x(a)
\phi_\pi(da|x)m_{\gamma,n}^{\phi_\pi}(dx\times
A)=\int_{S_1\setminus \hat{S}_1 \times A}
\frac{c_i(x,a)}{q_x(a)}M_{\gamma}^{n,\phi_\pi}(dx\times da)\\
&=&\int_{S_1\setminus \hat{S}_1 \times A}
\frac{c_i(x,a)}{q_x(a)}\textbf{P}_\gamma^{\sigma^M_\pi}(X_n\in
dx,~A_{n+1}\in da),
\end{eqnarray*}
where the second equality holds because of (\ref{ZyChapter57}),
and the third equality is by the definitions of
$M_{\gamma}^{n,\phi_\pi}$ and $m_{\gamma,n}^{\phi_\pi},$ and the
last equality is by (\ref{ZyChapter55}). For the second term in
the summand in the last line of (\ref{ZyChapter56}), we have
\begin{eqnarray*}
\int_{\hat{S_1}\times A} c_i(x,a)m_{\gamma,n}^{\phi_\pi}(dx\times
da)=\int_{\hat{S_1}\times A}
\frac{c_i(x,a)}{q_x(a)}\textbf{P}_\gamma^{\sigma^M_\pi}(X_n\in
dx,~A_{n+1}\in da),
\end{eqnarray*}
where the equality holds because
\begin{eqnarray*}
m_{\gamma,n}^{\phi_\pi}(\hat{S}_1\times
A)=0=\textbf{P}_\gamma^{\sigma^M_\pi}(X_n\in \hat{S}_1)
\end{eqnarray*}
with the first equality being by Lemma \ref{ZyChapterLem6} (see
(\ref{ZyChapter58}) therein) applied to $\phi_\pi$, which is
feasible with a finite value for problem (\ref{ZyChapter1}) for it
outperforms the policy $\pi$ by Corollary
\ref{ZyChapterCorollary3}, and the second equality being valid by
(\ref{ZyChapter55}) and that $M_\gamma^{n,\phi_\pi}(dx\times
da)=q_x(a)m_{\gamma,n}^{\phi_\pi}(dx\times da).$ For the third
term in the summand in the last line of (\ref{ZyChapter56}),
\begin{eqnarray*}
&&\int_{S_2\times A} c_i(x,a)m_{\gamma,n}^{\phi_\pi}(dx\times
da)=\int_{S_2} c_i(x,\psi^\ast(x))m_{\gamma,n}^{\phi_\pi}(dx\times
A)=0\\
&=&\int_{S_2}
\frac{c_i(x,\psi^\ast(x))}{q_x(\psi^\ast(x))}\textbf{P}_\gamma^{\sigma^M_\pi}(X_n\in
dx)=\int_{S_2\times A}
\frac{c_i(x,a)}{q_x(a)}\textbf{P}_\gamma^{\sigma^M_\pi}(X_n\in
dx,~A_{n+1}\in da),
\end{eqnarray*}
where the first equality is by (\ref{ZyChapter59}); recall that
$S_2\subseteq \zeta^c,$ the second and third equalities are by
(\ref{ZyChapter35}), and the last equality is by
(\ref{ZyChapter61}) and (\ref{ZyChapter23}). Finally, for the last
term in the summand of (\ref{ZyChapter56}), it holds that
\begin{eqnarray*}
&&\int_{S_3\times A} c_i(x,a)m_{\gamma,n}^{\phi_\pi}(dx\times
da)=\int_{S_3\times A}
\frac{c_i(x,a)}{q_x(a)}q_x(a)m_{\gamma,n}^{\phi_\pi}(dx\times
da)\\
&=&\int_{S_3\times A}
\frac{c_i(x,a)}{q_x(a)}M_{\gamma}^{n,\phi_\pi}(dx\times
da)=\int_{S_3\times A}
\frac{c_i(x,a)}{q_x(a)}\textbf{P}_\gamma^{\sigma^M_\pi}(X_n\in
dx,~A_{n+1}\in da),
\end{eqnarray*}
where the first equality is by the definition of $S_3,$ and the
last equality is by (\ref{ZyChapter55}). Combining these
observations, we see from (\ref{ZyChapter56}) that for each
$i=0,1,\dots,N,$
\begin{eqnarray}\label{ZyChapter62}
&&E_\gamma^{\phi_\pi}\left[\int_0^\infty \int_A
c_i(\xi_t,a)\phi_\pi(da|\xi_t)dt\right]=\sum_{n=0}^\infty
\int_{S\times A}
\frac{c_i(x,a)}{q_x(a)}\textbf{P}_\gamma^{\sigma^M_\pi}(X_n\in
dx,~A_{n+1}\in
da)\nonumber\\
&=&\textbf{E}_\gamma^{\sigma^M_\pi}\left[\sum_{n=0}^\infty
\frac{c_i(X_n,A_{n+1})}{q_{X_n}(A_{n+1})}\right].
\end{eqnarray}
On the other hand, one can apply Theorem 3.3 of Dufour et al
\cite{Dufour:2012} for the existence of a stationary policy
$\sigma^S_\pi$ for the DTMDP model $\{S_\infty,A,p,\gamma\}$
satisfying that for each $i=0,1,\dots,N,$
\begin{eqnarray*}
\textbf{E}_\gamma^{\sigma^S_\pi}\left[\sum_{n=0}^\infty
\frac{c_i(X_n,A_{n+1})}{q_{X_n}(A_{n+1})}\right]\le
\textbf{E}_\gamma^{\sigma^M_\pi}\left[\sum_{n=0}^\infty
\frac{c_i(X_n,A_{n+1})}{q_{X_n}(A_{n+1})}\right].
\end{eqnarray*}
This and (\ref{ZyChapter62}) thus prove the statement.
$\hfill\Box$

\begin{lemma}\label{ZyChapterLem8}
Suppose Condition \ref{BookWC1} is satisfied. Consider a
stationary policy $\sigma^S$ for the DTMDP model $\{S_\infty,
A,p,\gamma\}$, which satisfies
\begin{eqnarray*}
\sigma^S(da|x)=\delta_{f^\ast(x)}(da),~\forall~x\in S_2,
\end{eqnarray*}
and is optimal and with a finite value for problem
(\ref{ZyChapter63}). Here the transition probability $p(dy|x,a)$
is given by (\ref{ZyChapter64}) and (\ref{ZyChapter65}). Then
there is a stationary policy $\pi^S$ for the CTMDP problem
(\ref{ZyChapter1}) satisfying for each $i=0,1,\dots,N,$
\begin{eqnarray*}
\textbf{E}_\gamma^{\sigma^S}\left[\sum_{n=0}^\infty
\frac{c_i(X_n,A_{n+1})}{q_{X_n}(A_{n+1})}\right]=E_\gamma^{\pi^S}\left[\int_0^\infty
\int_A c_i(\xi_t,a)\pi^S(da|\xi_t)dt\right].
\end{eqnarray*}
\end{lemma}

\par\noindent\textit{Proof.} Since $\sigma^S$ is feasible with a finite value for problem
(\ref{ZyChapter63}), it is easy to see that
$\sum_{n=0}^\infty\textbf{P}_\gamma^{\sigma^S}(X_n\in
\hat{S}_1)=0$ so that, if necessary, we can modify the definition
of the policy $\sigma^S$ by putting
\begin{eqnarray*}
\sigma^S(da|x)=\delta_{\Delta}(da),~\forall~x\in \hat{S}_1,
\end{eqnarray*}
with $\Delta\in A$ being an arbitrarily fixed point; the resulting
policy is still optimal with a finite value for problem
(\ref{ZyChapter63}) and with the same performance vector as of the
original policy.

Note also that $\sigma^S(B(x)|x)=0$ for each $x\in S_1\setminus
\hat{S}_1.$ For this reason, we can legitimately define the
following stationary policy $\pi^S$ for the CTMDP model;
\begin{eqnarray*}
\pi^S(da|x)=\frac{\frac{1}{q_x(a)}\sigma^S(da|x)}{\int_A
\frac{1}{q_x(a)}\sigma^S(da|x)}
\end{eqnarray*}
for each $x\in S\setminus (S_2\bigcup \hat{S}_1),$
\begin{eqnarray*}
\sigma^S(da|x)=\delta_{\Delta}(da),
\end{eqnarray*}
for each $x\in \hat{S}_1,$ and
\begin{eqnarray*}
\pi^S(da|x)=\delta_{f^\ast(x)}(da)
\end{eqnarray*}
for each $x\in S_2.$ The discrete-time Markov chain $\{X_n\}$
under $\textbf{P}_\gamma^{\sigma^S}$ can be regarded as the
embedded chain of the pure jump time-homogeneous Markov process
$\{\xi_t\}$ under $P_\gamma^{\pi^S}$; see
\cite{FeinbergShiryayev:2014}. Indeed, it holds on ${\cal B}(S)$
that, for each $x\in S\setminus S_2,$
\begin{eqnarray*}
\int_A p(dy|x,a)\sigma^S(da|x)=\int_A
\frac{\tilde{q}(dy|x,a)}{q_x(a)}\sigma(da|x)=\frac{\int_A
\tilde{q}(dy|x,a)\pi^S(da|x)}{\int_A q_x(a)\pi^S(da|x)};
\end{eqnarray*}
for each $x\in S_2,$
\begin{eqnarray*}
\int_A p(dy|x,a)\sigma^S(da|x)=\int_A
\frac{\tilde{q}(dy|x,a)}{q_x(a)}\sigma^S(da|x)=\frac{\tilde{q}(dy|x,f^\ast(x))}{q_x(f^\ast(x))}=\frac{\int_A
\tilde{q}(dy|x,a)\pi^S(da|x)}{\int_A q_x(a)\pi^S(da|x)}=0;
\end{eqnarray*}
and for each $x\in \hat{S}_1,$
\begin{eqnarray*}
\int_A p(dy|x,a)\sigma^S(da|x)=\int_A
\frac{\tilde{q}(dy|x,a)}{q_x(a)}\sigma^S(da|x)=\frac{\tilde{q}(dy|x,\Delta)}{q_x(\Delta)}=\frac{\int_A
\tilde{q}(dy|x,a)\pi^S(da|x)}{\int_A q_x(a)\pi^S(da|x)}.
\end{eqnarray*}
Furthermore, it is easy to verify that for each $i=0,1,\dots,N,$
given the current state $x\in S$, the (conditional) expected total
cost during the current sojourn time of $\xi_t$ under
$P_\gamma^{\pi^S}$ is given by
$\int_A\frac{c_i(x,a)}{q_x(a)}\sigma^S(da|x)$, which is the same
as the (conditional) expected one-step cost for the discrete-time
Markov chain $\{X_n\}$ under $\textbf{P}_\gamma^{\sigma^S}.$ The
statement of this lemma now follows. $\hfill\Box$

\begin{condition}\label{FinitenessCon}
For problem (\ref{ZyChapter1}), there exists a feasible policy
with a finite value.
\end{condition}

\begin{theorem}\label{ZyChapterTheorem5}
Suppose Condition \ref{BookWC1} and Condition \ref{FinitenessCon}
are satisfied. Then for the CTMDP problem (\ref{ZyChapter1}),
there is a stationary optimal policy $\pi.$
\end{theorem}
\par\noindent\textit{Proof.} It is clear that for the CTMDP problem
(\ref{ZyChapter1}), one can be restricted to the class of feasible
policies $\pi$ with a finite value and satisfying
(\ref{ZyChapter53}); there exists at least one such policy under
Condition \ref{FinitenessCon}. It also holds that for the DTMDP
problem (\ref{ZyChapter63}), if the stationary policy $\sigma^S_1$
for the DTMDP model $\{S_\infty,A,p,\gamma\}$ is optimal, then the
stationary policy $\sigma^S$ for the DTMDP model
$\{S_\infty,A,p,\gamma\}$ defined by
$\sigma^S(da|x)=\sigma_1^S(da|x)$ for each $x\in S\setminus S_2$,
and $\sigma^S(da|x)=\delta_{f^\ast(x)}(da)$ for each $x\in S_2$ is
also optimal with a finite value for problem (\ref{ZyChapter63}).
Now the statement is a consequence of Lemma \ref{ZyChapterLem8},
Corollary \ref{ZyChapterCorollary4}, and Theorem 4.1 of Dufour et
al \cite{Dufour:2012}. $\hfill\Box$

\begin{remark}\label{ZyChapterRemark2}
Suppose Condition \ref{BookWC1} and Condition \ref{FinitenessCon}
are satisfied. Theorem 4.1 of Dufour et al \cite{Dufour:2012},
Lemma \ref{ZyChapterLem8} and Corollary \ref{ZyChapterCorollary4}
justify the reduction of problem (\ref{ZyChapter1}) for the CTMDP
model $\{S,A,q,\gamma\}$ to problem (\ref{ZyChapter63}) for the
DTMDP model $\{S_\infty,A,p,\gamma\}$; once the stationary optimal
policy for the DTMDP problem (\ref{ZyChapter63}), which exists, is
obtained, an optimal stationary policy for the CTMDP problem
(\ref{ZyChapter1}) can be automatically constructed based on it in
principle, and the two problems have the same value.
\end{remark}

\begin{remark}
As was rightly noted in \cite{Feinberg:2004}, if the transition
rates $q_x(a)$ are separated from zero, then one can show that for
each policy $\pi$ for the CTMDP, there is a policy $\sigma$ for
the DTMDP $\{S_\infty,A,p,\gamma\}$ such that
\begin{eqnarray*}
E_\gamma^\pi\left[\int_0^\infty \int_A
c_i(\xi_t,a)\pi(da|\omega,t)dt\right]=\textbf{E}_\gamma^\sigma\left[\sum_{n=0}^\infty
\frac{c_i(X_n,A_{n+1})}{q_{X_n}(A_{n+1})}\right]
\end{eqnarray*}
and vice versa, for each $i=0,1,\dots,N.$   The argument is
essentially the same as for the discounted case, and the reduction
is possible without further conditions. However, the objective of
the present paper is to consider the more delicate and nontrivial
case, i.e., when the transition rates are not necessarily
separated from zero.
\end{remark}

\section{Conclusion}\label{GuoZhangSec7}
To sum up, for the constrained total undiscounted optimal control
problem for a CTMDP in Borel state and action spaces, under the
compactness and continuity conditions, we showed the existence of
an optimal stationary policy out of the class of general
nonstationary ones. In the process, we justified the reduction of
the CTMDP model to a DTMDP model. Several properties about the
occupancy and occupation measures were obtained, too.

\appendix

\section{Appendix}

\subsection{Description of the DTMDP}
For the reader's convenience, we briefly describe the construction
of a DTMDP; see
\cite{Altman:1999,Bertsekas:1978,Hernandez-Lerma:1996,Piunovskiy:1997}
for greater details.

Consider a DTMDP with the (nonempty) Borel state space $X$,
(nonempty) Borel action space $U$, and the transition probability
$Q(dy|x,a)$, a stochastic kernel from $X\times U$ to ${\cal
B}(U).$ Let the initial distribution be $\mu(dx)$ on ${\cal
B}(X).$

A policy $\sigma$ for the DTMDP model $\{X,U,Q,\gamma\}$ is a
sequence of stochastic kernels $(\sigma_n)$, where for each
$n=1,2,\dots,$ $\sigma_n(du|z_0,u_1,\dots,u_{n-1},z_{n-1})$ is a
stochastic kernel from $X\times (U\times X)^{n-1}$ to ${\cal
B}(A)$. A policy $\sigma$ is called Markov if (with slight abuse
of notations) for each $n=1,2,\dots,$
$\sigma_n(du|z_0,u_1,\dots,u_{n-1},z_{n-1})=\sigma_n(du|z_{n-1}).$
A policy $\sigma$ is called stationary if (with slight abuse of
notations) for each $n=1,2,\dots,$
$\sigma_n(du|z_0,u_1,\dots,u_{n-1},z_{n-1})=\sigma(du|z_{n-1}).$ A
policy is called deterministic if all the stochastic kernels
$\sigma_n$ degenerate; if the stochastic kernels $\sigma_n$ do not
all degenerate, the policy is called randomized. A deterministic
stationary policy $\sigma$ with $\sigma(du|z)=\delta_{f(z)}(du)$,
where $f$ is a measurable mapping from $X$ to $U$, is often
denoted as $f$.

According to the Ionescu-Tulcea theorem,  under a policy
$\sigma=(\sigma_n)$, its strategic measure $\textbf{P}_\mu^\sigma$
is a probability measure on the countable product space $X\times
(U\times X)^\infty$ equipped with the Borel $\sigma$-algebra
defined by for each $\Gamma_U\in {\cal B}(U)$ and $\Gamma_X\in
{\cal B}(X)$,
\begin{eqnarray*}
&&\textbf{P}_\mu^\sigma(z_0\in \Gamma_X)=\mu(\Gamma_X),
\end{eqnarray*}
and for each $n=1,2,\dots,$
\begin{eqnarray*}
&&\textbf{P}_\mu^\sigma(u_n\in \Gamma_U|
z_0,u_1,z_1,\dots,u_{n-1},z_{n-1})=\sigma_n(\Gamma_U|z_0,u_1,z_1,\dots,u_{n-1},z_{n-1}),\\
&&\textbf{P}_\mu^\sigma(z_n\in \Gamma_X|
z_0,u_1,z_1,\dots,u_{n-1},z_{n-1},u_n)=Q(\Gamma_X|z_{n-1},u_n).
\end{eqnarray*}
The controlled and controlling processes are $\{z_n,u_{n+1}\}.$

\subsection{Auxiliary statements}
\begin{lemma}\label{ForwickChapterLem2}
Suppose Condition \ref{BookWC1}(b,c) is satisfied. Then the
following assertions hold.

\par\noindent(a) For each lower semicontinuous function $c(x,a)\in[0,\infty]$ on $S\times A,$ $\frac{c(x,a)}{q_x(a)}$ is lower semicontinuous in $(x,a)\in
S\times A;$  and $\int_S f(y)q(dy|x,a)$ is continuous in $(x,a)\in
S\times A$ for each bounded continuous function $f$ on $S.$

\par\noindent(b) For each $f(x)\in[0,\infty]$ lower
semicontinuous on $S,$ $\int_S f(y)\tilde{q}(dy|x,a)\in[0,\infty]$
is lower semicontinuous in $S\times A.$
\end{lemma}

\par\noindent\textit{Proof.} This statement is a consequence of
Lemma 7.14 of Bertsekas and Shreve \cite{Bertsekas:1978}, which
asserts that defined on a metric space, an extended real valued
lower semicontinuous function bounded from below is lower
semicontinuous if and only if there is an increasing sequence of
bounded continuous functions converging to it pointwise. We use
this fact without special reference below. Also recall the
convention of $\frac{0}{0}:=0$, $0\cdot\infty:=0,$ and
$\frac{x}{0}:=+\infty$ if $x>0$ (or $\frac{x}{0}:=-\infty$ if
$x<0$).

(a) Since $c(x,a)\in [0,\infty]$ is lower semicontinuous, there
exists a nondecreasing sequence of bounded continuous functions
$g_n(x,a)\ge 0$ on $K$ such that $\lim_{n\rightarrow
\infty}g_n(x,a)=c(x,a)$ for each $(x,a)\in S\times A.$  Since
$q_x(a)\ge 0,$ we see that $\frac{g_n(x,a)}{q_x(a)+\frac{1}{n}}\ge
0$ and increases to $\frac{c(x,a)}{q_x(a)}$ as $n\rightarrow
\infty$ for each $(x,a)\in S\times A.$ Remember that if
$c(x,a)=0$, then $g_n(x,a)=0$ for each $n$. Furthermore, it is
easy to see that
\begin{eqnarray*}
0\le \frac{g_n(x,a)}{ q_x(a)+\frac{1}{n}}\le n \sup_{(x,a)\in
K}\{g_n(x,a)\} <\infty
\end{eqnarray*}
and that $\frac{g_n(x,a)}{q_x(a)+\frac{1}{n}}$ is continuous in
$(x,a)\in S\times A$ for each $n$. Thus, $\frac{c(x,a)}{q_x(a)}$
is lower semicontinuous in $(x,a)\in S\times A.$

For the second assertion, we consider some bounded continuous
function $f$ on $S$. It follows from Condition \ref{BookWC1}(b,c)
that $\int_S f(y)q(dy|x,a)$ is continuous at each $(x,a)\in
S\times A$ such that $q_{x}(a)\ne 0.$ Now consider some
arbitrarily fixed $(x,a)\in S\times A$, where $q_x(a)=0,$ so that
$\frac{\int_S f(y)\tilde{q}(dy|x,a)}{q_{x}(a)} =0.$ Let some
convergent sequence $S\times A\ni(x_n,a_n)\rightarrow (x,a)\in
S\times A$ and some finite constant $\lambda>0$ be arbitrarily
fixed. If $\int_S f(y)\tilde{q}(dy|x_n,a_n)$ does not converge to
zero, then it would contradict
\begin{eqnarray*}
0=\lim_{n\rightarrow \infty}\frac{\int_S
f(y)\tilde{q}(dy|x_n,a_n)}{q_{x_n}(a_n)}= \lim_{n\rightarrow
\infty}\frac{\int_S
f(y)\tilde{q}(dy|x_n,a_n)}{\lambda+q_{x_n}(a_n)}\frac{\lambda+q_{x_n}(a_n)}{q_{x_n}(a_n)},
\end{eqnarray*}
where for the first inequality one should refer to Condition
\ref{BookWC1}(b), and for the contradiction we recall the
convention that $0\cdot\infty:=0,$ and note that
$\lim_{n\rightarrow
\infty}\frac{\lambda+q_{x_n}(a_n)}{q_{x_n}(a_n)}=+\infty.$ Thus,
$\int_S f(y)\tilde{q}(dy|x,a)$ is continuous at every $(x,a)\in
S\times A$. It follows from this and Condition \ref{BookWC1}(c)
that $\int_S f(y)q(dy|x,a)$ is continuous in $(x,a)\in S\times A$.
Since the bounded continuous function $f$ on $S$ is arbitrarily
fixed, this completes the proof of this part.

(b) Let $\{f^{(m)},~m=1,2,\dots\}$ be an increasing sequence of
nonnegative bounded continuous functions on $S$ such that
$f^{(m)}(x)\uparrow f(x)$ for each $x\in S.$ Then for each
$m=1,2,\dots,$
\begin{eqnarray*}
m\wedge\int_S f^{(m)}(y)\tilde{q}(dy|x,a)= m\wedge\left(\int_S
f^{(m)}(y)q(dy|x,a)- f^{(m)}(x)q_x(a)\right)\in [0,m]
\end{eqnarray*}
is a nonnegative bounded continuous function in $(x,a)\in S\times
A$ by the second assertion of part (a). Now the statement follows
from that $\{m\wedge\int_S
f^{(m)}(y)\tilde{q}(dy|x,a),~m=1,2,\dots\}$ is an increasing
sequence of nonnegative bounded continuous functions on $S$ such
that
\begin{eqnarray*}
\lim_{m\rightarrow \infty} m\wedge\int_S
f^{(m)}(y)\tilde{q}(dy|x,a)=\int_S f(y)\tilde{q}(dy|x,a),
\end{eqnarray*}
which is by the monotone convergence theorem. $\hfill\Box$

\begin{lemma}\label{ZyChapterLem3}
Suppose Condition \ref{BookWC1}(a) is satisfied, and let an
extended real-valued lower semicontinuous function $g$ on $S\times
A$ be fixed. Then the following assertions hold.
\par\noindent(a)  For each $\epsilon\in \mathbb{R},$ it holds that
the set $\left\{x\in S:~ \forall~a\in A,~g(x,a)>\epsilon\right\}$
is open in $S$.
\par\noindent(b)
\begin{eqnarray*}
\left\{x\in S: ~\forall~a\in A, ~
g(x,a)>0\right\}=\bigcup_{l=1}^\infty \left\{x\in S:~\forall~a\in
A,~g(x,a)>\frac{1}{l}\right\}.
\end{eqnarray*}
\end{lemma}
\par\noindent\textit{Proof.} See Lemmas 3.1 and 3.2 in
\cite{Dufour:2012}. $\hfill\Box$

\begin{lemma}\label{ZyChapterLem4}
Suppose Condition \ref{BookWC1}(a) is satisfied. Let
$\{g_n,~n=1,2,\dots\}$ be an increasing sequence of extended
real-valued functions on $S\times A$ such that for each
$n=1,2,\dots,$ $g_n(x,\cdot)$ is lower semicontinuous functions on
$A$ for each $x\in S.$ Then for each $x\in S,$
\begin{eqnarray}\label{ZyChapter8}
\lim_{n\rightarrow\infty}\inf_{a\in A}g_n(x,a)=\inf_{a\in
A}\lim_{n\rightarrow \infty}g_n(x,a).
\end{eqnarray}
\end{lemma}

\par\noindent\textit{Proof.}  See Proposition 10.1 of Sch\"{a}l \cite{Schal:1975a}; see also Appendix A of \cite{BauerleRieder:2011}.  $\hfill\Box$

\subsection{Some known facts about the unconstrained CTMDP problem}
Consider the following optimal control problem for the CTMDP model
$\{S,A,q\}$:
\begin{eqnarray}\label{ForwickChapter19}
E_x^\pi\left[\int_0^\infty \int_A \sum_{i=0}^N
c_i(\xi_t,a)\pi(da|\omega,t)dt\right]\rightarrow \min_{\pi\in
\Pi_{DM}},
\end{eqnarray}
the value function of which is denoted as
\begin{eqnarray}\label{ForwickChapter23}
V(x):=\inf_{\pi\in \Pi_{DM}} E_x^\pi\left[\int_0^\infty \int_A
c_0(\xi_t,a)\pi(da|\omega,t)dt\right].
\end{eqnarray}
Here $\Pi_{DM}$ stands for the class of deterministic Markov
policies for the CTMDP model $\{S,A,q\}.$ A policy $\pi^\ast\in
\Pi_{DM}$ is called optimal for problem (\ref{ForwickChapter19})
if
\begin{eqnarray*}
E_x^{\pi^\ast}\left[\int_0^\infty \int_A \sum_{i=0}^N
c_i(\xi_t,a)\pi^\ast(da|\omega,t)dt\right]=V(x)
\end{eqnarray*}
for each $x\in S.$

The following proposition is borrowed from \cite{Forwick:2004};
see Proposition 5.8 and Theorem 5.9 therein.
\begin{proposition}\label{ForwickChapterTheormeDP}
Suppose that Condition \ref{BookWC1} is satisfied. Then the
following assertions hold.
\par\noindent(a)
The function $ V$ is the minimal nonnegative lower semicontinuous
solution on $S$ to the following Bellman (optimality) equation:
\begin{eqnarray}\label{ForwickChapter16}
V(x)=\inf_{a\in A}\left\{\frac{\sum_{i=0}^Nc_i(x)}{q_x(a)}+\int_S
\frac{\tilde{q}(dy|x,a)}{q_x(a)}V(y)\right\}
\end{eqnarray}
for each $x\in S.$
\par\noindent(b) There is a deterministic stationary optimal policy $\varphi^\ast$
for the CTMDP problem (\ref{ForwickChapter19}), which can be taken
as a measurable mapping from $S$ to $A$ such that
\begin{eqnarray*}
\inf_{a\in A}\left\{\frac{\sum_{i=0}^Nc_i(x)}{q_x(a)}+\int_S
\frac{\tilde{q}(dy|x,a)}{q_x(a)} V(y)\right\}=
\frac{\sum_{i=0}^Nc_i(x,\varphi^\ast(x))}{q_x(\varphi^\ast(x))}+\int_S
\frac{\tilde{q}(dy|x,\varphi^\ast(x))}{q_x(\varphi^\ast(x))}V(y),~\forall~x\in
S.
\end{eqnarray*}
In fact, each deterministic stationary optimal policy for problem
(\ref{ForwickChapter19}) $\varphi^\ast$ satisfies the above
relation.
\end{proposition}
In fact, the authors of \cite{Forwick:2004} considered the more
general piecewise deterministic Markov decision process but in the
state space $\mathbb{R}^n$. When specializing to the case of a
CTMDP, one can put the more general Borel state space $S$.
Furthermore, the authors of \cite{Forwick:2004} assumed that
$V(x)<\infty$ for each $x\in S$; see ``Boundedness Assumption'' in
p.252 therein, which, could be withdrawn when specializing to the
CTMDP problem (\ref{ForwickChapter19}), as far as the validity of
the above proposition is concerned.

\begin{remark}\label{GuoZhangAppendixRem}
Proposition \ref{ForwickChapterTheormeDP} shows that under
Condition \ref{BookWC1} the CTMDP problem (\ref{ForwickChapter19})
can be reduced to the unconstrained DTMDP problem
\begin{eqnarray*}
\textbf{E}_x^\sigma\left[\sum_{n=0}^\infty \sum_{i=0}^N
\frac{c_i(X_n,A_{n+1})}{q_{X_n}(A_{n+1})}\right]\rightarrow
\min_{\sigma}
\end{eqnarray*}
because they have the same optimality equation. To this end, the
authors of \cite{Forwick:2004} firstly reduced the CTMDP model to
a more complicated DTMDP model, where the action space is in the
form of the space of all Borel measurable mappings, and had to
introduce the Young topology. The similar approach is developed in
\cite{BauerleRieder:2009,Costa:2013,Davis:1993,Yushkevich:1980}
for related unconstrained problems. On the other hand, the dynamic
programming approach based on the optimality equation is less
convenient for our constrained problem. In this paper, we choose
to study the occupation measures, and reduce directly the CTMDP
model to a DTMDP model with the same action space.
\end{remark}

\end{document}